\documentclass{tac}

\def\mathcaldef#1{\expandafter\def\csname#1\endcsname{{\cal#1}}}

\def\q{\quad}
\def\qq{\quad\quad}
\def\ot{\leftarrow}
\def\iff{\q\Leftrightarrow\q}
\def\open{\overleftarrow}
\def\closed{\overrightarrow}
\def\clopen{\overline}
\def\down{\downarrow\!\!}
\def\up{\uparrow\!\!}
\def\rup{\!\!\uparrow}
\def\rdown{\!\!\downarrow}
\def\bi{\updownarrow\!\!}
\def\rbi{\!\!\updownarrow}
\def\adj{\dashv}
\def\oX{\,^{\rm o}\! X}
\def\Xo{X^{\rm o}}
\def\oY{\,^{\rm o} Y}

\def\oXo{\,^{\rm o}\! X^{\rm o}}
\def\lp{\,^{\rm o}\pi}
\def\rp{\pi^{\rm o}}
\def\lpr{\,^{\rm o}\pi^{\rm o}}
\def\comp{\Gamma_{!}}
\def\disc{\Gamma^{*}}
\def\point{\Gamma_{*}}
\newcommand{\imp}{\mathop{\Rightarrow}\limits}
\newcommand{\iimp}{\mathop{\Longrightarrow}\limits}
\newcommand{\tto}{\mathop{\to}\limits}
\newcommand{\ttto}{\mathop{\longrightarrow}\limits}
\def\iso{\ttto^\sim}
\def\meets{\,\,\cap  \!\!\! !\,\,} 
\def\ch{\q\cdot\!\to\!\cdot\!\to\!\cdot\!\to\cdots\q}
\def\rch{\q\cdot\!\ot\!\cdot\!\ot\!\cdot\!\ot\cdots\q}
\def\fch{\q\cdot\!\to\!\cdot\!\to\cdots\to\!\cdot\q}
\def\dch{\q\cdots\to\!\cdot\!\to\!\cdot\!\to\cdots\q}
\def\op{^{\rm op}}
\def\pl{\pi_1}
\def\p2{\pi_2}
\def\phi{\varphi}
\def\eps{\varepsilon}

\mathcaldef{A}
\mathcaldef{B}
\mathcaldef{C}
\mathcaldef{D}
\mathcaldef{N}
\mathcaldef{Z}

\mathbfdef{Set}
\mathbfdef{Two}
\mathbfdef{Cat}
\mathbfdef{Pos}
\mathbfdef{Mod}
\mathbfdef{Bip}
\mathbfdef{T}
\mathbfdef{1}
\mathbfdef{2}
\mathbfdef{3}
\mathbfdef{Rel}

\def\Grf {{\rm\bf Grph}}
\def\Grfo {{\rm\bf Grph_0}}

\mathrmdef{hom}
\mathrmdef{ten}
\mathrmdef{true}
\mathrmdef{false}
\mathrmdef{P}
\mathrmdef{L}
\mathrmdef{ev}
\mathrmdef{Sp}
\mathrmdef{Ba}
\mathrmdef{At}
\mathrmdef{id}
\mathrmdef{gcd}
\mathrmdef{lcm}
\mathrmdef{tot}

\let\pf\proof
\let\epf\endproof

\newtheorem{prop}{Proposition}
\newtheorem{lemma}{Lemma}
\newtheorem{corol}{Corollary}

\author{Claudio Pisani}
\address{via Gioberti 86,\\ 10128 Torino, Italy.}
\title{Bipolar Spaces}
\copyrightyear{2005}
\keywords{discrete fibrations and opfibrations, reflections and coreflections, tensor, cofigures, negation, atoms, 
Cauchy completion, Kan extensions, mappings}
\amsclass{18Axx}
\eaddress{pisclau@yahoo.it}

\begin{document}

\maketitle

\begin{abstract}
Some basic features of the simultaneous inclusion of discrete fibrations and discrete opfibrations on a category $\A$
in the category of categories over $\A$ are studied; in particular, the reflections and the coreflections of the latter in 
the former are considered, along with a negation-complement operator which, applied to a discrete fibration, gives a functor
with values in discrete opfibrations (and vice versa) and which turns out to be classical, in that the strong contraposition law holds.
Such an analysis is developed in an appropriate conceptual frame that encompasses similar    
``bipolar" situations and in which a key role is played by ``cofigures", that is components of products; 
e.g. the classicity of the negation-complement operator corresponds to the fact that discrete opfibrations (or in general ``closed 
parts") are properly analyzed by cofigures with shape in discrete fibrations (``open parts"), that is, that the latter are ``coadequate" 
for the former, and vice versa. 
In this context, a very natural definition of ``atom" is proposed and it is shown that, in the above situation, the category of atoms
reflections is the Cauchy completion of $\A$.
\tableofcontents
\end{abstract}

\section{Introduction}
\label{introduction}

Every functor $F:\A\to\B$ gives rise to an interpretation of any object $B$ of $\B$ as a presheaf on $\A$ 
(the ``category of shapes"), whose $A$-elements, namely the arrows $f:FA\to B$ in $\hom_\B(FA,B)$, 
may be called figures of $B$ with shape $A\in\A$.

The main thesis of the present paper is that in some important situations that naturally arise in mathematics 
an explicit consideration of the functor 
\[ \ten_\B:\B\times \B\to\Set \] 
defined by components of products, is conceptually advantageous. 
It plays in some sense a complementary or dual role to the 
$\hom_\B$ functor and should be thought of as the truth value of the proposition ``$A$ meets $B$", in the same sense according to which
${\hom}_\B(A,B)$ may be tought of as the truth value of the proposition ``$A$ is included in $B$" (see Sections~\ref{two} and~\ref{preorders}).
Specifically, in the above situation, the functor $F$ gives rise to another interpretation of any object $B$ of $\B$ 
as a (covariant) presheaf on $\A$, the ``cointerpretation" of $B$ via $F$, whose $A$-elements, namely the components of $FA\times B$ in 
$\ten_\B(FA,B)$, may be called ``cofigures" of $B$ with shape $A$. If the cointerpretation functor is full and faithful we say that
$F$ is ``coadequate".

For example, if $\B$ is $\Grf$, the category of (irreflexive) graphs, then appropriate choices of $F$ lead to the interpretation of a graph
as various sets (such as the sets of its nodes, its arrows, its loops) or as various mappings between these sets (such as the inclusion of 
loops in arrows, or the domain and codomain mappings) and so on. On the other hand, a graph can be cointerpreted as the set of its components, as 
the set of its nodes (once again, since the dot graph is an ``atom"; see later), as the mapping which takes a node to its component, and so on. 
While it is easy to find ``small" adequate subcategories of $\Grf$ (typically the full one generated by the dot and the arrow), 
even the identity functor is not coadequate for $\Grf$ (e.g., a graph cannot be distinguished from its thin 
reflection by ``tensoring" with other graphs). On the other hand, the ``right endomappings", that is the graphs with a bijective domain mapping, 
{\em are} coadequate for the subcategory of ``left endomappings" (with bijective codomain). Furthermore, if $\N$ is the 
monoid of the natural numbers and $F:\N\to\Grf$ is the functor given by the action of $\N$ on the infinite chain $\ch$ by translations, 
then the corresponding interpretation and cointerpretation give the coreflection of graphs in right endomappings and the reflection in left endomappings 
respectively (this example is treated in detail in Section~\ref{endomappings}). 

Since various situations display analogous features, it seems appropriate to condense them in the following definition      
of a {\bf bipolar space} \[ X = \langle \P(X),\open X,\closed X\rangle \]
\begin{enumerate}
\item
$\P(X)$ is a category with tensor (see Section~\ref{components})
and $\open X$ and $\closed X$ are full replete subcategories of $\P(X)$.
The objects of $\P(X)$ are the {\bf parts} of the space $X$ and those in $\open X$ or in $\closed X$ are called {\bf open} and {\bf closed} parts, 
or also ``left mappings" and ``right mappings", respectively.
\item \label{neg}
If $A\in\open{X}$, then there is a {\bf complement} or {\bf negation} functor $\neg A$
right adjoint to $\ten(A,-)$: \[ \ten(A,-)\adj \neg A:\Set\to \P(X) \] 
which takes values in $\closed X$;
and dually, any closed part has an open complement (that is, with open values).
Often the bipolar spaces considered are {\bf strong}, in that the following stronger property holds: 
for any open $A$ and closed $D$, the exponential object $A\imp D$  exists in $\P(X)$ and is itself closed, and vice versa.  

\item
Open and closed parts are reciprocally coadequate.
As shown in Section~\ref{coadequacy}, this is equivalent to the fact that
the complement-negation is classical, in that the strong contraposition law
\[ \frac{\neg B \to \neg A}{A \to B} \]
holds, for $A$ and $B$ both open or both closed; so open or closed parts can be ``recovered" from their nagation.
 
\item
$\open X$ and $\closed X$ are both reflective and coreflective in $\P(X)$, with reflections 
\[ \down(-):\P(X)\to\open X \qq\qq \up(-):\P(X)\to\closed X \] 
and coreflections 
\[ (-)\rdown\, :\P(X)\to\open X \qq\qq (-)\rup\, :\P(X)\to\closed X \]
As a consequence we have the natural isomorphisms
\begin{eqnarray} \label{f1}
\hom(P,Q\rdown) \cong \hom(\down P,Q\rdown) \cong \hom(\down P,Q) \\
\hom(P,Q\rup) \cong \hom(\up P,Q\rup) \cong \hom(\up P,Q)    \nonumber
\end{eqnarray}
\begin{equation} \label{f2}
\ten(\down P,Q) \cong \ten(\down P,\up Q) \cong \ten(P,\up Q) 
\end{equation} 
for parts $P$ and $Q$ of $X$; while the former are obvious, the latter follow from 
point~\ref{neg} above, as will be shown in Section~\ref{components}.
      
\end {enumerate}

\noindent The last axiom is in fact a very strong requirement; e.g. it excludes, in the two-valued case, general topological spaces 
(see Section~\ref{preorders}). Anyhow, it should be clear which hypothesis are really needed in the various part of the subsequent theory.

The key example is $\Sp X = \langle \Cat/X,\open X,\closed X\rangle$, the bipolar space associated to a category $X$, 
which has $\P(\Sp X) = \Cat/X$, while $\open{\Sp X}=\open X$ and $\closed{\Sp X}=\closed X$ are the full subcategories of discrete fibrations and 
discrete opfibrations over $X$ respectively; so $\Sp X$ has categories over $X$ as parts, and a part is open (respectively closed) iff it is
a discrete fibration (respectively, a discrete opfibration). 

Among the parts of $\Sp X$, the objects $x:{\bf 1}\to X$ of $X$ have the special property that for any open part $A$ and any closed part $D$
\begin{equation} \label{f3}
\hom(x,A) \cong\ten(A,x) \qq \qq \hom(x,D) \cong\ten(x,D)
\end{equation}
since all of them reduce to the objects or components of a discrete category, namely the fibre of $A$ or $D$ over $x$; 
moreover, both the bijections are mediated by the same universal element in $\ten(x,x)$ and so any object $x\in X$ is an ``atom" 
(see Sections~\ref{atoms} and~\ref{atomcat}).
The reflections $\down x$ and $\up x$ are the discrete fibrations $X/x$ and $x/X$ of objects over and under $x$, and correspond to the 
representable functors $X(-,x)$ and $X(x,-)$, the required adjunctions being provided by the Yoneda lemma.
Similarly, any idempotent $e$ in $X$, as a part of $\Sp X$, is an atom whose reflections $\down e$ and $\up e$ are the adjoint bimodules 
associated to $e$ in the Cauchy completion of $X$, and we shall see in Section~\ref{atomcat2} that the same is true for the reflections 
of any atom $x\in\At(\Sp X)$.

Similarly, in a two-valued context, for any poset $X$ we have the bipolar space $\Sp X = \langle {\cal P}(X),\open X,\closed X\rangle$  
which has as parts the parts (that is, subsets) of the underlying set of $X$, and a part is open (respectively closed) iff it is a sieve 
(respectively, a cosieve). Here the (strong) atoms are the parts with a single element.

If $X$ is a graph, we also have a bipolar space $\Sp X = \langle \Grf/X,\open X,\closed X\rangle$, having as open (respectively closed) parts
those graphs over $X$ which are ``left" (respectively ``right") functional over any arrow of $X$.
Here the atoms are the nodes of $X$ as graphs over $X$.

In general, given a bipolar space $X$, there are two categories $\oX$ and $\Xo$ with the atoms $x\in\At(X)$ as objects and with  
$\oX(x,y) = \hom(\down x,\down y)$ and $\Xo(x,y) = \hom(\up x,\up y)$, and the obvious full and faithful 
functors $\lp:\oX\to\P(X)$ and $\rp:\Xo\to\P(X)$. Furthermore, there is a duality 
\begin{equation} \label{du}
\sigma:\Xo\iso (\oX)\op 
\end{equation}
and the corresponding interpretation and cointerpretation of open parts as presheaves on $\oX$ and as covariant presheaves on $\Xo$ are 
equivalent modulo this isomorphism. Indeed, for any open part $A$ and any atom $x$, equations~(\ref{f1}),~(\ref{f2}) 
and~(\ref{f3}) above induce bijections
\[ \hom(\down x,A) \cong\ten(A,\up x) \]
that define a natural isomorphism
\[ \hom(\lp(\sigma -),A)\iso\ten(A,\rp -):\Xo\to\Set \] 
as shown in Section~\ref{atoms}.
Symmetrically, any closed part $D$ is (co)interpreted as a functor $\oX\to\Set$ in two distinct but equivalent ways:
\[ Dx:=\hom(\up x,D) \qq\qq Dx:=\ten(\down x,D) \]
If the bipolar space $X$ is ``atomic", that is if $\oX$ is adequate for $\open X$ and $\Xo$ is adequate for $\closed X$ 
(or, equivalently, $\oX$ is coadequate for $\closed X$ and $\Xo$ is coadequate for $\open X$), then open and closed parts ``are" 
presheaves on $\oX$ and $\Xo$.
In particular, the reflection $\up P$ of a part $P\in\P(X)$ ``is" the functor $\oX\to\Set$ defined by
\[ (\up P) x := \ten(\down x,\up P) \qq\qq (\up P)\alpha := \ten(\alpha,\up P) \] 
for any $\alpha:\,\down x\to\,\down y$ in $\oX$, and similarly $\down P$ ``is" a functor $\Xo\to\Set$:
\[ (\down P) x := \ten(\down P,\up x) \qq\qq (\down P)\beta := \ten(\down P,\beta) \] 
for any $\beta:\,\up x\to\,\up y$ in $\Xo$. Using the~(\ref{f2}) we get the isomorphisms of functors
\begin{eqnarray} \label{r1}
(\up P) x \cong \ten(\down x,P) \qq\qq (\up P) \alpha \cong \ten(\alpha,P) \\
(\down P) x \cong \ten(P,\up x) \qq\qq (\down P) \beta \cong \ten(P,\beta) \nonumber
\end{eqnarray} 
Similarly, using the interpretation instead of the cointerpretation and the~(\ref{f1}) instead of the~(\ref{f2}), we get the following 
formulas for the coreflections:
\begin{eqnarray} \label{r2}
(P\rdown) x \cong \hom(\down x,P) \qq\qq (P\rdown) \alpha \cong \hom(\alpha,P) \\ 
(P\rup) x \cong \hom(\up x,P) \qq\qq (P\rup) \beta \cong \hom(\beta,P) \nonumber
\end{eqnarray}
So, in an atomic space $X$, $\up P$ is obtained as the cointerpretation of $P$ via $\lp:\oX\to\P(X)$, while $P\rup$ is obtained as the interpretation
of $P$ via $\rp:\Xo\to\P(X)$ (and dually for $\down P$ and $P\rdown$).

This applies to $\Sp X$, where $X$ is a category or a graph. In the case of graphs we find in particular the formulas for the reflection
and the coreflection of graphs in endomappings, touched upon at the beginning of this introduction. Indeed, if $\L$ is the terminal graph  
the only atom in $\Sp\L$ is the dot ${\rm D}\in\Grf$, its reflections are the chains $\ch$ or $\rch$ and $\oX$ and $\Xo$ are the full subcategories 
of their endomorphisms, which as categories are isomorphic to the monoid $\N$ of natural numbers.

\subsubsection*{Outline}
After introducing categories with tensor in Section~\ref{tensor}, where also the two-valued case is considered, and atoms in 
Section~\ref{atoms}, in the next three sections we discuss examples of bipolar spaces.

Sections~\ref{graphs} and~\ref{preorders} are mainly intended as a sort of extended introduction. In the former
we present two particularly pregnant instances of spaces $\Sp X$ associated to a graph $X$, namely we consider $\Sp{\rm A}$ 
and $\Sp\L$, where ${\rm A},\L\in\Grf$ are the arrow and the loop respectively.
The latter is devoted to the most typical kind of bipolar spaces in the two-valued context, the spaces $\Sp X$ associated to a poset $X$;
they coincide essentially with the Alexandrov topological spaces.

The central example, which may be seen as the Alexandrov space associated to a category $X$, is presented in Section~\ref{categories},
where in particular the properties of atoms in $\Sp X$ are proved.
We show in Section~\ref{reflections} that the formulas~(\ref{r1}) and~(\ref{r2}) above give indeed the expected 
reflections and coreflections of a category $P$ over $X$ in the discrete (op)fibrations.

In Section~\ref{doctrines} we define ``continuous maps", the arrows of the category $\Bip$ of bipolar spaces, so that $\Sp X$ 
becomes the object mapping of a (pseudo)functor $\Sp:\Cat\to \Bip$. In the other direction, $\oX$ (or $\Xo$) will also be extended to a 
functor $\Ba:\Bip\to \Cat$. As sketched before, $\Ba(\Sp X)$ gives ``the" Cauchy completion of the category $X$; coherently,
in the two-valued context, $\Ba(\Sp X)$ is equivalent to the poset $X$.
Similarly, there are functors $\Sp:\Grf\to \Bip$ and $\Sp:\Grfo\to \Bip$, where $\Grfo$ denotes the category of reflexive graphs;
in these cases $\Ba(\Sp X)$ gives the free category on $X$.
The connections with Kan extensions of set functors will be also dealt with. 

In the Appendix, after briefly considering reflections and coreflections in discrete bifibrations (and more generally in ``clopen" parts), 
we outline to what extent a definition of functionality via cofigures differs from the usual one via orthogonality of figures, and finally
we hint at possible directions for further research.

\subsection*{Related works}
The study of similar ``bipolar" situations, and in particular of that of Section~\ref{categories}, has been carried out from a different
viewpoint in~\cite{bun00}; the fact that discrete fibrations on a category $\A$ are reflective in the category of categories over $\A$
was treated in~\cite{par73} and~\cite{law70} (see also~\cite{law73}); the fact that they are coreflective too, 
although presumably widely known as well, has not appeared in a published form to the author's knowledge.
The general attitude is in the spirit of Lawvere's ``generalized logic" (see~\cite{law73} and~\cite{law86}); also the viewpoint of 
Section~\ref{doctrines} is cleary borrowed from Lawvere (see in particular~\cite{law70}), whose ideas have strongly influenced those of 
the author; the resulting doctrines $\Sp X$ are a sort of temporal doctrines (see Section~\ref{abstractions}).

\section{Categories with tensor}
\label{tensor}

In this section and in the next one we build the conceptual frame sketched in the introduction, by developing the necessary tools.

\subsection{Components, tensor and negation} 
\label{components}

We say that a category $\A$ {\bf has components} (or that $\A$ is a category {\bf with components}) if there are functors
\[ \comp \adj \disc \adj \point :\A\to\Set \]
called the {\bf components}, the {\bf discrete} and the {\bf points} functor, respectively.
The objects of $\A$ isomorphic to those of form $\disc S$ (for a set $S\in\Set$) are said to be {\bf discrete} or {\bf constant}, 
while elements of $\point A$ are called points and those of $\comp A$ {\bf components} (or also {\bf copoints}) of the object $A\in\A$.
If $\A$ is a category with components, then it is such in an essentially unique way, since in this case
$\disc 1$ is terminal in $\A$, $\disc S$ is the copower $S\cdot 1$ (that is the sum of $S$ copies of $1\in\A$) 
and the points functor is represented by $1\in\A$.
 
If $\A$ is a category with components which is also cartesian (meaning that it has finite limits), 
we define the functor 
\[ \ten_\A : \A\times\A\to\Set \] 
by components following products, that is
\[ \ten_\A(A,B) = \comp (A\times B) \]
and we say that $\A$ {\bf has tensor} (or that $\A$ is a category {\bf with tensor}).
The elements of $\ten_\A(A,B)$ can be called {\bf cofigures} of $B$ with shape $A$ (or also, because of the symmetry of $\ten$,
cofigures of $A$ with shape $B$). Note that we are using the prefix ``co" also to distinguish the concepts 
which refer to the $\ten$ functor from the corresponding ones relative to the $\hom$ functor.  

Any presheaf category has tensor, $\point$  and $\comp$ being the limit and colimit functors to $\Set$;
in particular, if $M$ is a monoid and $X$ is an $M$-set, the points of $X$ are the fixed points of the corresponding action of $M$
and its components are the orbits. If $M$ is commutative, then the unit $\eta_X:X\to\disc\comp X$
has a specially simple form: $\eta_X x = \eta_X y$ iff $mx = ny$ ($m,n\in M$), while, if $M$ is a group, $\eta$ has the even simpler
form $\eta_X x = \eta_X y$ iff $y = gx$ ($g\in M$). 
As an exemple, consider a group $G$ and the left action of $G$ on itself; as any representable functor, 
the $G$-set $G$ is connected (that is $\comp G = 1$) while $\ten(G,G)$ is (the underlying set of) $G$.
More interestingly, if $N$ is the monoid of natural numbers with addition, $\ten(N,N)$ 
can be seen as the set of integer numbers since $\eta(m,n) = \eta(m',n')$ iff either $m+h=m'$ and $n+h=n'$ 
or $m'+h=m$ and $n'+h=n$, for a natural number $h$. Similarly, if $N^*$ is the multiplicative monoid of non zero natural numbers
then $\ten(N^*,N^*)$ can be seen as the set of (non zero) rational numbers, obtained as a quotient of the 
set of fractions $N^*\times N^*$.

If $\A$ is also (cartesian) closed, that is if there is an adjunction 
\[ -\times A \adj A{\imp}- \]
with parameter $A\in\A$,
then the dual or complementary roles played by the functors $\ten_\A$ and $\hom_\A$ 
are more evident, as summarized in the following  
\begin{prop} \label{closed}
Let $\A$ be a closed category with tensor, then
\begin{eqnarray}
\label{1}
\ten(A,B) := \comp (A\times B) \q &&\q \ten_\A(1,-) \cong \comp \\
\label{2}
\hom(A,B) \cong \point (A\imp B) \q && \q \hom_\A(1,-) \cong \point  \\ 
\label{3} 
\ten(-,A) &\adj & A\imp \disc(-) \\ 
\label{4}
\disc(-)\times A &\adj & {\hom}(A,-) 
\end{eqnarray}
the last two being adjunctions with parameter $A\in\A$.
\end{prop}
\epf
Note that equations~(\ref{1}) and~(\ref{2}) hold in the non closed case too, with the exception of the first of~(\ref{2}) 
which may have no meaning. In general, the copower functor $-\cdot A:\Set\to\A$, left adjoint to the representable 
$\hom(A,-)$, may exist without being given by the formula displayed in equation~(\ref{4}) (as in the case of categories
with a zero object, when the functors $\comp$, $\disc$ and $\point$ become constant). On the other hand, if the right adjoint 
to the ``corepresentable functor"
\[ \ten(-,A) \adj \neg A : \Set\to\A \]
exists, then it is easily seen to have the form~(\ref{3}). In this case, that is if $A\in\A$
is exponentiable for constant bases, we call the functor $\neg A$ the {\bf negation} or {\bf complement}
of $A$, and say that $A$ ``has negation" or ``complement" (in~\cite{law96} a similar definition of negation is given). 
Note that $\neg 1 = \disc$. 

It is worth stressing that non isomorphic objects may corepresent isomorphic functors (and so, in particular, may have
isomorphic negations), since a category need not to be coadequate for itself (see Section~\ref{coadequacy}). 
An example, already mentioned in the Introduction, is given by two graphs in $\Grf$ with the same thin reflection or, 
more dramatically, by any two {\em reflexive} graph with the same number of components (see Section~\ref{space}).
In the two-valued setting (see Section~\ref{two}) more familiar examples are given by elements of Heyting algebras with the 
same pseudocomplement (such as the open subsets of ${\bf\rm R}$ obtained taking off a finite number of points).

In presence of distributivity, the property of having tensor becomes stable:
\begin{prop} \label{slice}
Let $\A$ be a distributive category with tensor, then any slice $\A/A$ has tensor. Furthermore, the tensor
functor $\ten_\A : \A\times\A\to\Set$ preserves sums in each argument: 
\[ \ten(A+B,C) \cong \ten(A,C)+\ten(B,C) \qq\q \ten(A,B+C) \cong \ten(A,B)+\ten(A,C) \]
\end{prop}
\pf Since $\A/A$ is still cartesian and has copowers, we only have to show that there is a
component functor $\comp$ left adjoint to $\disc$. Indeed observe that, by distributivity, $\disc S\times A$ gives the 
copower $S\cdot 1$ of the terminal $\id_A$ in $\A/A$; so, by composition of adjoints, $\comp^{\A/A}=\comp^\A\circ\tot$, 
where $\tot:\A/A\to\A$ is the ``total" functor, left adjoint to $A\times -:\A\to\A/A$. 
The proof of the last statement is immediate, since both $A\times -$ and $\comp$ preserve sums. 
\epf

\subsection{Tensor and negation in bipolar spaces} 
\label{bip}

We now illustrate some of the interactions between the axioms for bipolar spaces, in particular those 
concerning the tensor functor and (co)reflectivity.
We begin with a proof of the ``coadjunction" property mentioned in the introduction (Corollary~\ref{coadj}) 
and conclude with some remarks about exponentiability. 
\begin{lemma}
Let $\A$ be a category with tensor and let $\A'$ be a full reflective subcategory of $\A$, with reflection $(-)':\A\to\A'$; 
if $B\in\A$ has a negation-complement $\neg B$ with values in $\A'$, then the unit $\eta$ of the reflection induces (natural) bijections
\[ \ten(\eta_A,B) : \ten(A,B) \iso \ten(A',B) \]
\end{lemma}
\pf The result follows from Yoneda and the following chain of bijections, natural in all variables:
\[ \Set(\ten(A,B),S) \cong \A(A,(\neg B)S) \cong \A(A',(\neg B)S) \cong \Set(\ten(A',B),S) \]
A more careful inspection shows that the bijections are in fact induced by the unit of the reflection.
\epf
Since in a bipolar space $X = \langle \P(X),\open X,\closed X\rangle$
every open part in $A\in\open X$ has a closed complement $\neg A: \Set\to \closed X$, and conversely,
and since the full subcategories $\open X$ and $\closed X$ of $\P(X)$ have reflections   
$\down (-):\P(X)\to\open X$ and $\up (-):\P(X)\to\closed X$, then the above lemma applies to give
\begin{prop} \label{pcoadj}
Given a bipolar space $X$, there are bijections 
\[ \ten(A,P) \cong \ten(A,\up P) \qq\qq \ten(P,D) \cong \ten(\down P,D) \]
natural in $P\in\P(X)$, $A\in\open X$ and $D\in\closed X$. 
Furthermore, the bijections are induced by the units $\lambda_P$ and $\rho_P$ of the open and closed 
reflections $\down(-)$ and $\up(-)$ respectively. 
\end{prop} \epf
\begin{corol} \label{coadj}
If $P$ and $Q$ are parts of a bipolar space, there are natural isomorphisms
\[ \ten(\down P,Q) \cong \ten(\down P,\up Q) \cong \ten(P,\up Q) \]
\end{corol} \epf
\noindent So $\up(-)$ and $\down(-)$ are ``coadjoint" functors. Note however that the remark following Proposition~\ref{closed} implies that
coadjoint functors do not determine each other up to isomorphisms.

In Proposition~\ref{slice} we showed a stability property of categories with tensor;
another one is the following
\begin{prop} \label{proptop}
If $\A$ is a category with tensor and $\A'$ is a reflective and coreflective full replete subcategory of $\A$, then $\A'$ has 
itself tensor, the points, discrete and components functors being the ``same" as those of $\A$. 
\end{prop}
\pf Indeed, $\A'$ is closed with respect to the limits and colimits which exist in $\A$ and so these serve as limits and colimits for $\A'$ 
too (see~\cite{bor94}). In particular, $\A'$ is cartesian and has copowers, and so the discrete functor $\disc:\Set\to\A$ takes values in $\A$; 
the points and the components functor for $\A'$ can then be obtained by restricting those of $\A$.
\epf
\noindent Then in the above situation the inclusion functor preserves the tensor functor (as well as the $\hom$ functor): 
if $A$ and $B$ are objects of $\A'$  
\[ \hom_{\A'}(A,B) = \hom_\A(A,B) \qq\qq \ten_{\A'}(A,B) = \ten_\A(A,B) \]
Returning to bipolar spaces, the above proposition applies to the subcategories $\open X$ and $\closed X$ of $\P(X)$, giving
\begin{corol} \label{top}
The categories of open and closed parts of a bipolar space $X$ have themselves tensor, and the
points, components and tensor functors therein are the restriction of the corresponding ones in $\P(X)$. 
Limits and colimits of open or closed parts which exist in $\P(X)$ are again open or closed, respectively;
in particular the constant-discrete parts are both open and closed. 
\end{corol} \epf
In the same hypothesis of Proposition~\ref{proptop}, the inclusion of $\A'$ in $\A$ may not 
preserve the exponentials that exist in $\A'$; actually, if $A$ and $B$ are objects of $\A'$ such that $A\imp^\A B$ 
exists, then its coreflection in $\A'$ gives the exponential $A\imp^{\A'} B$. So we have
\begin{prop}   \label{clopen}
In a bipolar space $X = \langle \P(X),\open X,\closed X \rangle$, 
if $B$ and $D$ are closed parts and $B\imp D$ exists in $\P(X)$, then $(B\imp D)\rup$ gives the exponential in $\closed X$.
In particular, if $X$ is a strong bipolar space, $B$ is a clopen part (that is, both closed and open) and $D$ is closed, 
then the closed part $B\imp D$ is also the exponential $B\imp D$ in $\closed X$. 
Symmetric properties hold for open parts.
\end{prop}
\epf
\begin{corol}   \label{clopencor}
The full subcategory $\clopen X$ of clopen parts of a strong bipolar space $X$ is a cartesian closed category with tensor, 
and the inclusion of $\clopen X$ in $\P(X)$ preserves the tensor functor and exponentials.
\end{corol}
\epf

\subsection{Adequacy and coadequacy}  
\label{coadequacy}

We now define coadequacy, which plays with respect to the tensor functor the same role that adequacy
plays with respect to the $\hom$ functor, and we will show that, as mentioned in the introduction, the reciprocal coadequacy
between open and closed parts of a bipolar space corresponds to the classicity of negation therein (Corollary~\ref{contr2}).

Let $T:\A\times\B\to\C$ a functor, to be thought of as a generalized bimodule. We say that $T$ is {\bf left adequate} if the
corresponding $\A\to [\B,\C]$ is full and faithful (where $[\B,\C] := \B\iimp^\Cat\C$).
Similarly, if $\B\to [\A,\C]$ is full and faithful, then $T$ is {\bf right adequate}.
More generally, if in the above situation a functor $F:\B'\to\B$ and a full subcategory $\A'$ of $\A$ are also given, 
we say that $F$ is left adequate for $\A'$ (via $T$) if the corresponding $\A'\to [\B',\C]$ is full and faithful 
(and similarly for the dual right notion). 
If $T = \hom_\A:\A\op\times\A\to\Set$ we find the usual notion of Isbell adequacy,
while if $T = \ten_\A:\A\times\A\to\Set$ we have the corresponding notion of {\bf coadequacy},
which is basic in the study of bipolar spaces. (Warning: some authors use the term ``coadequate" meaning ``right adequate").
Now we turn to the relations between coadequacy and negation.
\begin{lemma}
If $L:\A\times\B\to\C$ and $R:\B\op\times\C\to\A$ are functors such that 
\[ L(-,B) \adj R(B,-) \]
is an adjunction with parameter $B\in\B$, that is if there are bijections
\[ \C(L(A,B);C)\iso \A(A;R(B,C)) \]
natural in $A$, $B$ and $C$, then there are natural bijections 
\[ [\A,\C](L(-,B);L(-,B'))\iso [\C,\A](R(B',-);R(B,-)) \]
\end{lemma}
\pf We have the following chain of (natural) bijections:
\begin{eqnarray*}
[\A,\C ](L(-,B);L(-,B'))\cong \int_{A\in\A} \C(L(A,B);L(A,B')) \\
\cong \int_{A\in\A} \int_{C\in\C} \Set[\C(L(A,B'),C);\C(L(A,B),C) ] \\
\cong \int_{C\in\C} \int_{A\in\A} \Set[\A(A;R(B',C));\A(A;R(B,C)) ] \\
\cong \int_{C\in\C} \A(R(B',C);R(B,C)) \cong [\C,\A ](R(B',-);R(B,-)) 
\end{eqnarray*}
where the end formula for natural transformations, the interchange property of ends and the Yoneda lemma have been used
(see~\cite{mac71}). 
\epf
\begin{prop} \label{contr}
Given an adjunction with parameter $L(-,B) \adj R(B,-)$ as above, $L$ is right adequate iff $R$ is left adequate. 
\end{prop}
\pf Indeed, by the above lemma, both conditions reduce to 
\[ [\A,\C ](L(-,B);L(-,B'))\cong \B(B,B') \cong [\C,\A ](R(B',-);R(B,-)) \]
\epf
Since every open part of a bipolar space $X = \langle \P(X),\open X,\closed X \rangle$ has a closed complement, and conversely,
we have the following adjunctions with parameters $A\in\open X$ and $D\in\closed X$ respectively (see Proposition~\ref{closed}):
\[ \ten(-,A) \adj \neg A : \Set \to \closed X \qq\qq \ten(-,D) \adj \neg D : \Set \to \open X \] 
and since furthermore open and closed parts are assumed to be reciprocally coadequate, Proposition~\ref{contr} above gives
\begin{corol} \label{contr2}
In a bipolar space, the strong contraposition law 
\[ \frac{\neg B \to \neg A}{A \to B} \]
holds, for $A$ and $B$ both open or both closed; that is, there is a (natural) bijection between the natural tansformations
above the horizontal line and the inclusions (morphisms) of parts below it. 
\end{corol} 
\epf

\subsection{The two-valued case}  
\label{two}

The parallel analysis of a (categorical) concept in the two-valued and in the set-valued contexts has the advantage
of a mutual clarification: the former gives a semplified perspective and a better intuitive grasp; 
the latter allows a deeper understanding, since it displays the ``solid" concept and not only its (often blurred) shadow.
We now specialize the above concepts to the two-valued setting, replacing the category of truth values $\Set$ with $\Two$,
the poset (isomorphic to the arrow category $\2$) which has ``$\true$" and ``$\false$" as objects and $\false\leq\true$.

Here, by a poset we mean a $\Two$-valued category $\A$ (also known as preordered set).
As usual, we interpret $\A(A,B)$ as the truth value of the proposition $A\leq B$, or ``$A$ is included in $B$".
We denote by $\top$ (respectively $\bot$) a terminal (initial) object of a poset $\A$, that is a maximum (minimum) of $\A$.
The poset $\A$ has components  
\[ \comp \adj \disc \adj \point :\A\to\Two \] 
if and only if it has both a maximum and a minimum. Indeed, if this is the case we have

\begin{enumerate}
\item
$\disc (\true) = \top\q$ and $\q\disc (\false) = \bot$
\item
$\point A = \true \iff \top\leq A$
\item
$\comp A = \false \iff A\leq\bot$
\end{enumerate}
So we may interpret the points functor $\point$ as giving the truth value of the proposition ``$A$ is a maximum of $\A$" or ``$A$ is full",
and the components functor $\comp$ as giving the truth value of the proposition ``$A$ is not a minimum of $\A$" or ``$A$ is not empty".

A poset $\A$ with tensor is then simply a meet semilattice that has a minimum. 
If this is the case, we denote $\ten(A,B)$ by $A\meets B$: 
\[ A\meets B := \comp (A\wedge B) \] 
and intrepret it as the truth value of the proposition ``$A$ meets $B$" (that is, their meet is not empty). 
The following proposition follows directly from the above definitions, or from the corresponding $\Set$-valued facts in Subsection~\ref{components}.
\begin{prop} \label{pseudo}
An object $A$ of a poset with tensor $\A$ has negation $\neg A:\Two\to\A$ iff it has a pseudocomplement $A\imp\bot$.
In this case, $\neg A (\false)$ is the pseudocompement of $A$ itself.
\end{prop}
\epf
In the two-valued context, a bipolar space $X = \langle \P(X),\open X,\closed X\rangle$ consists of a poset $\P(X)$ with tensor 
(that is, a meet semilattice with $\bot$) and two sub(po)sets $\open X$ and $\closed X$
(of open and closed parts) such that 
\begin{enumerate}
\item
every open part has a closed (pseudo)complement, and conversely; 
\item
the strong contraposition law for open and for closed parts hold; this is equivalent to the double negation law
\[ \neg\neg A\leq A \] 
for $A$ open or closed, and, as in the set-valued case (and as can be easily proved directly), also 
to the fact that the open parts are coadequate for the closed ones and conversely: 
\[ \frac{\forall D\in\closed X \quad A\meets D\q\imp\q B\meets D}{A \leq B} \qq\qq
     \frac{\forall A\in\open X \quad A\meets D\q\imp\q A\meets E}{D \leq E} \] 
\item
$\open X$ and $\closed X$ and both reflective and coreflective in $\P(X)$.
\end{enumerate}
The space $X$ is strong if for any open part $A$ and closed part $D$ the relative pseudocomplements $A\imp D$ and $D\imp A$ exist 
and are respectively closed and open.

The two-valued version of Corollary~\ref{coadj} is the following simple formula:
\[ \down P\meets\, Q \iff \down P\meets\up Q \iff P\meets\up Q \]
for any parts $P,Q\in \P(X)$. 

\subsection{Discrete spaces and Boolean categories}  
\label{boole}

Say that a bipolar space $X = \langle \P(X),\open X,\closed X\rangle$ is {\bf Boolean} if open and closed parts coincide, 
that is $\open X =\closed X$ as subcategories of $\P(X)$; among Boolean spaces there are two extremal cases, with a maximum 
or a minimum of clopen (that is, both open and closed) parts: the {\bf discrete} spaces, in which any part is clopen, 
and the {\bf codiscrete} spaces, in which the only clopen parts are the constant-discrete ones. The discrete spaces have the form 
\[ X = \langle \P(X),\P(X),\P(X)\rangle \] 
that is, they ``are" categories $\P(X)$ with tensor which are coadequate for themselves and such that any object is exponentiable for 
discrete objects; $X$ is strong when the corresponding category $\P(X)$ is cartesian closed.

In the two-valued case the discrete spaces are the pseudocomple\-mented semilattices with classical pseudocomplement.
Strong discrete bipolar spaces are then the classical Brouwerian logics (in the sense of~\cite{eik66}) and so they are
(equivalent to) the Boolean algebras.

The above considerations point toward a possible determination of what a ``Boolean category" should be, 
at least to the extent that the proportion ``posets are to categories what Boolean algebras are to Boolean categories" 
is supposed to hold: it seems reasonable to identify Boolean categories with the strong discrete bipolar spaces,
that is, with cartesian closed categories with components that are coadequate for themselves; indeed the posetal ones
are the Boolean algebras, as just remarked.
The typical Boolean categories are the presheaf categories $\open X$ on a grupoid $X$;
indeed, as shown in Section~\ref{categories}, if $X$ is a category, then $\Sp X = \langle \Cat/X,\open X,\closed X\rangle$ is a 
strong bipolar space, and so, if $X$ is a grupoid, then $\open X=\closed X$ is coadequate for itself (see Corollary~\ref{clopencor}).
Conversely, it seems likely that the category of presheaves on $X$ is Boolean only if $X$ is a grupoid, but the standard proof of the 
corresponding fact in the two-valued context cannot be staightforwardly extended, resting on the law of excluded middle which is not
available in the set-valued case.

\section{Atoms in bipolar spaces}
\label{atoms}

In this section, after recalling some very basic categorical concepts such as universal elements of bifunctors, we present
atoms and their relevant features. In Section~\ref{atomcat2} we shall prove that, in the bipolar space $\Sp X$
associated to a category $X$, the reflections of atoms are the retracts of representable functors, while in Section~\ref{base}
the good functorial behaviour of atoms with respect to continuous maps of bipolar spaces is analyzed.

\subsection{Universal elements of bifunctors}
\label{biuniversal}

Let $t:\A\times\B \to \Set$ be a set-valued bifunctor and $u\in t(A,B)$; we say that $u$ is universal for $A$ if 
it is a universal element of the functor $t(A,-):\B\to \Set$ that is if $\beta\mapsto t(A,\beta)u$ gives bijections
\[ t(A,-)u : \hom_\B(B,B') \iso t(A,B') \] 
for any $B'\in\B$, so defining a natural isomorphism $\hom_\B(B,-) \iso t(A,-)$.
Given an object $A\in A$, we say that $A$ has a universal element if there is an object $B\in\B$ and an element $u\in t(A,B)$ universal
for $A$; this is of course equivalent to the fact that $t(A,-)$ is a representable functor. 
Similarly one defines universal elements for $B$ and says that $B$ has a universal element.
If $u\in t(A,B)$ is universal for both $A$ and $B$, it is called {\bf biuniversal}, and we say that $A$ or $B$ have a biuniversal element.
The following are well-known facts:
\begin{enumerate}
\item
Given two universal elements of $A$, $u\in t(A,B)$ and $u'\in t(A,B')$, there is a unique isomorphism $\beta:B\iso B'$ such that $t(A,\beta)u = u'$.
\item
If $\A'$ is the full subcategory of the objects of $\A$ that have universal elements, then one can define a functor $F:\A'\to\B\op$ 
such that there is a natural isomorphism between the restriction $t:\A'\times \B \to \Set$ and $\hom_\B(F-,-)$.
Namely, $FA$ is any of the objects of $\B$ such that there is $u_A\in t(A,FA)$ universal for $A$, and once such a family $(u_A, A\in\A')$
is fixed, $F$ is uniquely defined on arrows $\alpha:A\to A'$ in $\A'$ by 
\[ t(A',F\alpha)u_{A'} = t(\alpha,FA)u_A \in t(A',FA) \]
\item  \label{dual}
If $(u_i\in t(A_i,B_i),i\in I)$ is a family of biuniversal elements, there is a duality between the full subcategories $\A_0$ and $\B_0$
generated by the $A_i$ and the $B_i$; the duality becomes an isomorphism $\A_0\cong\B_0\op$, rather than just an equivalence, if 
$\A_0$ and $\B_0$ are taken as the categories with objects in $I$ and arrows in $\A$ and $\B$ respectively.
\end{enumerate}
Familiar consequences of these remarks are easily drawn by dualizing one of the arguments of $t$,
e.g. considering a bifunctor $t:\A\op\times\B\to\Set\q$: 
\begin{itemize}
\item
if any object $A\in\A$ has a universal element, then $t$ ``is" a functor $\A\to\B$ (defined up to isomorphism);
\item
if furthermore any object $B\in\B$ also has a universal element, then $t$ ``is" an adjunction $\A\rightharpoonup\B$;
\item
the adjuction restricts to an equivalence between the full subcategories $\A_0$ and $\B_0$ generated by the objects of $A\in\A$ and $B\in\B$
which have biuniversal element, or equivalently, by the objects of $A\in\A$ such that the unit $\eta_A$ of the adjunction is an isomorphism 
and by the objects $B\in\B$ such that the counit $\eps_B$ is an isomorphism. 
\item
in particular, if any universal element is in fact biuniversal, then $t$ ``is" an equivalence $\A\simeq \B$.
\end{itemize}

\subsection{Atoms and atomic pairs}
\label{atomic}

The idea behind our definition of atom is very simple: in the two-valued context, a strong atom of a bipolar space is a part $x\in\P(X)$
``so small" that it meets any part $P$ iff it is included in $P$ itself, that is
\[ x\meets P \iff x\leq P \]
(but not ``too small": the minimum $\bot$ is included in every part but doesn't meet any part, since $\bot\meets P \Leftrightarrow\comp\bot = \false$).
In the set-valued context, the condition becomes: $x$ meets any part $P$ in as many ways as it is included in $P$ itself, that is
\[ \ten(x,P) \cong \hom(x,P) \]
More precisely, a part $x\in\P(X)$ of a bipolar space $X = \langle \P(X),\open X,\closed X\rangle$ is a {\bf strong atom} if there is a biuniversal element 
\[ u\in\ten(x,x) \]
for the bifunctor $\ten:\P(X)\times\P(X)\to\Set$. Then the above mentioned bijection is mediated by $u$:
\[ \ten(x,-)u : \hom(x,P) \iso \ten(x,P) \]
for any $P\in\P(X)$. This is in fact a very strong requirement, and often the biuniversal property holds for open 
and closed parts only: a part of $x$ of $X$ is an {\bf atom}  if there is $u\in\ten(x,x)$ such that $\alpha\mapsto\ten(\alpha,x)u$ is a bijection 
\[ \ten(-,x)u : \hom(x,A) \iso \ten(A,x) \]
for any $A\in\open X$ and also $\beta\mapsto\ten(x,\beta)u$ is a bijection 
\[ \ten(x,-)u : \hom(x,D) \iso \ten(x,D) \] 
for any $D\in\closed X$.
Slighty improperly, we say that $u\in\ten(x,x)$ is biuniversal for the atom $x$ also in this weaker case. 
$\At X$ is the set of the atoms of the space $X$.

A pair $\langle\down x,\up x\rangle \in \open X\times\closed X$ is called an {\bf atomic pair} of $X$ if there is a biuniversal element 
\[ u\in\ten(\down x,\up x) \] 
for the restriction of the tensor functor $\ten:\open X\times\closed X\to\Set$.
As the notation suggests, the ``left" (that is open) and ``right" (closed) reflections of an atom form an atomic pair, the {\bf reflection} of $x$:
\begin{prop}  \label{atomic1}
If $x\in\P(X)$ is an atom of $X$, then the pair $\langle\down x,\up x\rangle$ is atomic.
\end{prop}
\pf Let $u\in\ten(x,x)$ be biuniversal for $x$ and define 
\[ u' = \ten(\lambda_x,\rho_x)u\in\ten(\down x,\up x) \]
where $\lambda$ and $\rho$ denote the units of the left and right reflections $\down(-)$ and $\up(-)$.
To prove that $u'$ is the desired biuniversal element, let us check e.g. the right universality: the mapping 
\[ \ten(\down x,-)u' : \hom(\up x,D) \to \ten(\down x,D) \]
takes $\beta\in\hom(\up x,D)$ to 
\[ \ten(\down x,\beta)u' = \ten(\down x,\beta)\ten(\lambda_x,\rho_x)u = \ten(\lambda_x,\beta\circ\rho_x)u \]
and then factorizes through
\[ \beta \mapsto \beta\circ\rho_x \mapsto \ten(x,\beta\circ\rho_x)u \mapsto \ten(\lambda_x,\beta\circ\rho_x)u \]
where the first is a bijection 
\[ \hom(\rho_x,D) : \hom(\up x,D) \iso \hom(x,D) \] 
by the universality of the unit $\rho_x$, the second is a bijection 
\[ \ten(x,-)u : \hom(x,D) \iso \ten(x,D) \] 
by the (bi)universality of $u$, while the third is a bijection 
\[ \ten(\lambda_x,x) : \ten(x,D)\iso \ten(\down x,D) \] 
by the ``universality" of $\lambda_x$ with respect to $\ten$, following from Proposition~\ref{pcoadj}.
\epf
As sketched in the Introduction, there are categories $\oX$ and $\Xo$ with objects in $\At X$ and arrows in $\P(X)$, 
according to the restrictions 
\begin{equation}   \label{atomiceq2}
\down(-)\, ,\, \up(-) : \At x\to\P(X)
\end{equation}
and ``projections" functors 
\[ \lp:\oX\to\P(X) \qq\qq \rp:\Xo\to\P(X) \]
with the~(\ref{atomiceq2}) as object mappings.
Alternatively, one could have defined the categories $\oX$ and $\Xo$ as the full subcategories
of $\P(X)$ generated by the left and right reflections of atoms; the resulting theory would be essentially the same, but with the drawback of
having an equivalence instead of the following isomorphism (see remark~(\ref{dual}) of Section~\ref{biuniversal}):
\begin{prop}  \label{atomic2}
The categories $\oX$ and $\Xo$ are dual.
\end{prop}
\pf Indeed, selecting a biuniversal element $u_x$ for any atom $x$ of $X$, there are bijections 
\[ \hom(\down x,\down y) \iso \ten(\down y,\up x) \qq\qq  \hom(\up y,\up x) \iso \ten(\down y,\up x) \]
given by
\[ \alpha\mapsto\ten(\alpha,\up x)u'_x \qq\qq \beta\mapsto\ten(\down y,\beta)u'_y \]
where $u'_x$ and $u'_y$ are the corresponding universal elements of the reflection $\langle\down x,\up x\rangle$ of $x$, 
as in Proposition~\ref{atomic1}. Then by composition we have bijections 
\[ \sigma_{x,y} : \hom(\down x,\down y) \iso \hom(\up y,\up x) \]
wherein two morphism $\alpha:\down x\to\down y$ and $\beta:\up y\to\up x$ in $\P(X)$ correspond to each other iff 
\[ \ten(\alpha,\up x)u'_x = \ten(\down y,\beta)u'_y \in \ten(\down y,\up x) \]
To prove that this defines an isomorphism of categories 
\[ \sigma : \Xo \iso (\oX)\op \]
it suffices to observe that the identity corresponds to itself, and that if $\alpha:\down x\to\down y$ 
corresponds to $\beta:\up y\to\up x$ and $\alpha':\down y\to\down z$ corresponds to $\beta':\up z\to\up y$
in above bijections, then $\alpha'\circ\alpha:\down x\to\down z$ corresponds to $\beta\circ\beta':\up z\to\up x$ because
\begin{eqnarray*} 
\ten(\alpha'\circ\alpha,\up x)u'_x = \ten(\alpha',\up x)\circ\ten(\alpha,\up x)u'_x \\
= \ten(\alpha',\up x)\circ\ten(\down x,\beta)u'_y = \ten(\alpha',\beta)u'_y = \ten(\down x,\beta)\circ\ten(\alpha',\up x)u'_y \\
= \ten(\down x,\beta)\circ\ten(\down x,\beta')u'_z = \ten(\down x,\beta\circ\beta')u'_z 
\end{eqnarray*}
\epf
\begin{prop}    \label{atomic3}
The interpretation of an open part via $\lp:\oX\to\P(X)$ and its cointerpretation via $\rp:\Xo\to\P(X)$ are equivalent
according to the above isomorphism $\sigma:\Xo \iso (\oX)\op$; symmetrically, the interpretation of a closed part via $\rp$ 
and its cointerpretation via $\lp$ are equivalent according to $\sigma$.
\end{prop}
\pf Considering the case of closed parts, we want to show that the bijections
\[ \hom(\up x,D) \iso \ten(\down x,D) \]
mediated by the biuniversal element $u'_x\in\ten(\down x,\up x)$ define a natural isomorphism
\[ \hom(\rp(\sigma^{-1} -),D) \iso \ten(\lp -,D):\oX\to\Set \] 
Indeed, if $\alpha:\down x\to\down y$ corresponds to $\beta:\up y\to\up x$ under the isomorhism $\sigma$ and $\xi\in\hom(\up x,D)$, 
then following the two paths from $\hom(\up x,D)$ to $\ten(\down y,D)$ in the square 
\[ \xi \mapsto \hom(\beta,D)\xi = \xi\circ\beta \mapsto \ten(\down y,\xi\circ\beta)u'_y \] 
\[ \xi \mapsto \ten(\down x,\xi)u'_x \mapsto \ten(\alpha,D)(\ten(\down x,\xi)u'_x) \]
we get the same result, because 
\begin{eqnarray*} 
\ten(\alpha,D)(\ten(\down x,\xi)u'_x) = \ten(\alpha,\xi)u'_x = \ten(\down y,\xi)(\ten(\alpha,\up x)u'_x) \\
= \ten(\down y,\xi)(\ten(\down y,\beta)u'_y) = \ten(\down y,\xi\circ\beta)u'_y 
\end{eqnarray*}
\epf
\begin{corol}   \label{atomic4}
$\lp:\oX\to\P(X)$ is adequate for open parts iff $\,\rp:\Xo\to\P(X)$ is coadequate for open parts. 
Symmetrically, $\rp:\Xo\to\P(X)$ is adequate for closed parts iff $\lp:\oX\to\P(X)$ is coadequate for closed parts.
\end{corol}
If both the conditions above hold, we say that the bipolar space $X$ is {\bf atomic}.
So, in an atomic bipolar space the open parts $A$ ``are" presheaves on $\oX$ (or covariant presheaves on $\Xo$), 
according to any of the following formulas
\begin{eqnarray} 
A x := \hom(\down x,A) &;& A \alpha := \hom(\alpha,A):\hom(\down y,A)\to\hom(\down x,A) \\
A x := \ten(A,\up x) &;& A \beta := \ten(A,\beta):\ten(A,\up y)\to\ten(A,\up x) 
\end{eqnarray} 
with $\alpha:\down x\to\down y$ corresponding to $\beta:\up y\to\up x$ under the isomorphism $\sigma : \Xo \iso (\oX)\op$;
symmetrically, the closed parts $D$ ``are" presheaves on $\Xo$ (or covariant presheaves on $\oX$) in two distinct but equivalent ways:
\begin{eqnarray} 
D x := \hom(\up x,D) &;& D \beta := \hom(\beta,D):\hom(\up x,D)\to\hom(\up y,D) \\
D x := \ten(\down x,D) &;& D \alpha := \ten(\alpha,D):\ten(\down x,D)\to\ten(\down y,D) 
\end{eqnarray}
As explained in the Introduction, Corollary~\ref{coadj} then gives  
\begin{corol}    \label{atomic5}
If $X$ is atomic, the cointerpretation of a part $P\in\P(X)$ via $\lp:\oX\to\P(X)$ gives the reflection $\up P$ in $\closed X$, while
the interpretation of $P$ via $\rp:\Xo\to\P(X)$ gives the coreflection $P\rup$ in $\closed X$, and symmetrically for the reflection $\down P$ 
and the coreflection $P\rdown$ of $P$ in open parts.
\end{corol}
\epf
\noindent To be more explicit, the formulas~(\ref{r1}) and~(\ref{r2}) in the Introduction hold, where the $\alpha$ and $\beta$ are arrows 
in $\oX$ and $\Xo$ respectively, and so ultimately they are all inclusions of parts in $\P(X)$. 

The above discussion readily specializes to the two-valued setting of Section~\ref{two}. In particular we have the dual 
posets $\oX$ and $\Xo$ associated to a two-valued bipolar space $X$, and the related concept of atomic space;
in such a space, any open (respectively closed) part ``is" a sieve (respectively cosieve; see Section~\ref{preorders}) on $\oX$,
and the (co)reflections of a part $P$ can be obtained by the atoms reflections and the application of the $\hom$ and $\ten$ 
functors, which now reduce to the ``is included" and ``meets" predicates (see Section~\ref{two}):
\begin{equation} \label{refp}
x\in\,\up P \iff \down x\meets\, P  \qq\qq  x\in\,\down P \iff P\meets\up x 
\end{equation}
and similarly for the coreflections:
\begin{equation}  \label{corp}
x\in P\rdown \iff \down x\leq P  \qq\qq  x\in P\rup \iff \up x\leq P
\end{equation}

\subsection{Evaluation at atoms}
\label{eval}

For any atom $x$ of a bipolar space $X$, we can consider the functor 
\[ \ev_x:\closed X\to\Set \] 
obtained by restricting the representable $\hom(x,-):\P(X)\to\Set$ or the corepresentable $\ten(x,-):\P(X)\to\Set$ 
to closed parts:
\[ \ev_x D := \hom(x,D) \cong \ten(x,D) \cong \hom(\up x,D) \cong \ten(\down x,D) \]
and dually for open parts. Thus $\ev_x$ is represented by $\up x$, and it is the restriction to closed parts of the functor 
\[ \ten(\down x,-):\P(X)\to\Set \]
corepresented by an open part, and so it has both a left adjoint and a right adjoint, given respectively by copowers and negation
(see Proposition~\ref{closed}): 
\[ -\cdot\up x:\Set\to\closed X \qq\qq  \neg\down x:\Set\to\closed X \] 
In particular, $\ev_x$ preserves limits and colimits and so $\up x$ is an absolutely presentable 
object of $\closed X$. The following examples may become clearer after reading some of the subsequent sections: 
\begin{enumerate}
\item
In the case of the bipolar space $\Sp\L$ associated to the terminal graph (see Section~\ref{endomappings}), 
the ``dot" atom $\rm D$ leads to the functor $\ev_{\rm D}$ which takes a (right) endomapping to
the set on which it acts. Then the above formulas give its left and right adjoint functors, which take a set $S$ to 
the copower $S\cdot {\rm C}_{\infty}$ of infinite chains, or to the set of sequences in $S$ respectively.
This example is easily generalized to $\Sp X$ for any graph $X$ (see Section~\ref{space}).
\item
If $X$ is a two-valued space, 
\[ \ev_x D := x\leq D \iff x\meets D \iff \up x\leq D \iff \down x\meets D \]
that is, the evaluation of $D$ at $x$ gives the truth value of the proposition ``$x$ is included in $D$"; 
its left adjoint takes ``$\true$" to $\up x$, while its right adjoint takes ``$\false$" the pseudocomplement of $\down x$
in $\P(X)$. Perhaps it is worth to see directly why $\up x$ is a (completely) prime element of $\closed X$: 
\begin{eqnarray*}
\up x \,\leq\, D\vee E \q &\vdash &\q \false \\
\down x \,\meets\, (D\vee E) \q &\vdash &\q \false \\
D\vee E \q &\leq &\q \neg\down x\,(\false) \\
D \q\leq\q \neg\down x\,(\false) \q &\& &\q E \q\leq\q \neg\down x\,(\false) \\
\down x \,\meets\, D \q\vdash\q \false \q &\& &\q \down x \,\meets\, E \q\vdash\q \false \\
\up x \,\leq\, D \q\vdash\q \false \q &\& &\q \up x \,\leq\, E \q\vdash\q \false \\
\up x \,\leq\, D \, &{\rm or} &\, \up x \,\leq\, E \q\vdash\q \false
\end{eqnarray*}
so that
\[ \up x \,\leq\, D\vee E \q\iff\q \up x \,\leq\, D \q {\rm or}\q \up x \,\leq\, D \] 
where $\vdash$, $\&$ and ``or" denote respectively the $\hom$ functor, the product and the coproduct in the category $\Two$, 
while $\leq$ and $\vee$ denote the $\hom$ functor and the coproduct in the $\Two$-category $\P(X)$; recall that $\neg D(\false)$
is the pseudocomplement of $D$ in $\P(X)$ (see Proposition~\ref{pseudo}). 
This is an example of proof in the ``generalized logic with component" (see the Conclusions), since it is valid in the set-valued,
or also $\cal V$-valued setting too: just replace ``false" with the generic truth value in $\cal V$ and use Yoneda in the last step.

In particular, if $X$ the two-valued discrete bipolar space corresponding to a Boolean algebra, as in Section~\ref{boole},
then an atom $x$ is a completely prime element of $X$, and so it is an atom in the usual sense too. 
\item
In the case of the two-valued bipolar space $\Sp X$ associated to a poset $X$ of Section~\ref{preorders},
we have that, for an atom $x\in X$, $\ev_x:\closed X\to\Two$ takes an upper set (or cosieve) $D$ to the truth value of 
the proposition ``$x\in D$"; its left adjoint takes ``$\true$" to the principal cosieve $\up x$, while its right adjoint
takes ``$\false$" to the complement in $\P(\Sp X) = {\cal P}(X)$ of $\down x$, that is, to the set of $y\in X$
such that $y\not\leq x$. 
\item
In the case of the bipolar space $\Sp X$ associated to a category $X$ of Section~\ref{categories}, as outlined in the Introduction
there are two kinds of atoms (actually, the former are a particular case of the latter): those corresponding to the objects of $X$ 
and those corresponding to the idempotent arrows in $X$ (see Sections~\ref{atomcat} and~\ref{atomcat2}).
In the former case, for any object $x\in X$ and (covariant) presheaf $D\in\closed X$, $\ev_x D$ is simply $Dx$, 
that is $\ev_x:\closed X\to\Set$ is the usual evaluation functor. Its left adjoint takes a set $S\in\Set$ to the functor 
\[ S\cdot\up x = S\times X(x,-) : X\to\Set \]
while its left adjoint $\neg\down x$ is given by
\[ (\neg\down x)S = \Set(X(-,x),S) : X\to\Set \]
On the other hand, for any idempotent arrow $e:x\to x$ in $X$, the evaluation $\ev_e D$ of a (covariant) presheaf $D$ at the corresponding 
atom gives the elements in $Dx$ fixed by $e$. The values of its left and right adjoint at a set $S\in\Set$ now are 
\[ S\cdot(\up e) = S\times X'(x,-) \qq\qq (\neg\down e)S = \Set(X'(-,x),S) \]
where $X'(x,y)$ is the set of arrows $f:x\to y$ such that $f\circ e = f$, and symmetrically for $X'(y,x)$ (see Section~\ref{atomcat2}).
\end{enumerate}

\section{Graphs and mappings}
\label{graphs}

This section is devoted to two examples of bipolar spaces $\Sp X$ associated to a graph $X$, namely those associated to
the arrow and to the loop graphs; although the general case will be considered in Section~\ref{doctrines}, we treat these special cases 
in some detail because they are noteworthy and paradigmatical in many respects.

The category $\Grf$, as any presheaf category, has tensor: given $P\in\Grf$, $\point P$ is the set of its loops, while $\comp P$ is the set 
of its components, given by the coequalizer of the domain and codomain mappings 
(if $\eta_P:P\to\disc\comp P$ is the unit of the adjunction, and $x$ and $y$ are nodes of $P$, then $\eta_P x = \eta_P y$ iff
there is an undirected path from $x$ to $y$).

It is convenient to fix some notations for $\Grf$. We denote by $\rm L$ the loop (with one node and one arrow), which is the terminal graph, 
by $\rm D$ the dot (with one node and no arrows) and by $\rm A$ the arrow (with two nodes and one arrow between them). 
The two morphisms from the dot to the arrow are indicated by ${\rm D}\tto^{\delta_0}{\rm A}$ and ${\rm D}\tto^{\delta_1}{\rm A}$.
${\rm C}_n$ is the chain $\fch$ with $n$ nodes (e.g. ${\rm C}_0$ is the empty graph, ${\rm C}_1$
is the dot and ${\rm C}_2$ is the arrow) and ${\rm C}_\infty$ is the infinite chain $\ch$; 
note that the ${\rm C}_n$'s are self-dual (that is, ${\rm C}_n\op\cong {\rm C}_n$) while ${\rm C}_{\infty}\op=\rch$.

\subsection{Graphs over the arrow and mappings}
\label{mappings}

Mappings can be seen as particular graphs over the arrow $\rm A$, that is objects of $\Grf/{\rm A}$.
In fact, we now define the bipolar space 
\[ \Sp{\rm A} = \langle \Grf/{\rm A},\open{\rm A},\closed{\rm A}\rangle \]
which has the graphs over $\rm A$ as parts, while $\open{\rm A}$ and $\closed{\rm A}$ are the full subcategories of ``left" and ``right" mappings,
and in so doing we show how the framework of Sections~\ref{tensor} and~\ref{atoms} give a new perspective on some familiar facts about mappings.

To keep on with the notational policy of using uppercase letters to denote parts, we (locally) define ${\rm A}:={\rm A}\tto^{\rm id_A}{\rm A}$,
${\rm D_0}:={\rm D}\tto^{\delta_0}{\rm A}$ and ${\rm D_1}:={\rm D}\tto^{\delta_1}{\rm A}$, while 
${\rm D_0}\tto^{\delta_0}{\rm A}$ and ${\rm D_1}\tto^{\delta_1}{\rm A}$ denote the (unique) morphisms in $\Grf/{\rm A}$ to the terminal $\rm A$.

An object $P = P\tto^{\pi}{\rm A}$ of $\Grf/{\rm A}$ can be described by the sets $P(0,1)$ (the arrows of $P$), $P(0)$ and $P(1)$
(the nodes of $P$ over each node of $A$) and the mappings domain and codomain $p_0:P(0,1)\to P(0)$ and $p_1:P(0,1)\to P(1)$.
A morphism $\phi : P\to Q$ is given by three mappings $\phi(0,1):P(0,1)\to Q(0,1)$, $\phi(0):P(0)\to Q(0)$ and
$\phi(1):P(1)\to Q(1)$ that commute with $p_0$ and $p_1$. (Alternatively, one can consider, instead of $P(0,1)$, a ``matrix" of 
sets $P(x,y)$ ($x\in P(0)$, $y\in P(1)$) like in Section~\ref{catexp}.)

We say that $P$ is a {\bf right mapping} ({\bf left mapping}) if $p_0$ (respectively $p_1$) is bijective, or, equivalently, if $P$ is interpreted as a bijection
by the arrow $\delta_0$ (respectively $\delta_1$) in $\Grf/{\rm A}$, that is $\delta_0$ (respectively $\delta_1$) is orthogonal to $P$.
So, if $P$ is a right mapping, from each element $x\in\P(0)$ ``starts" a unique arrow toward $P(1)$, and so it ``is" a mapping
$P:P(0)\to P(1)$, and similarly a left mapping $Q\in\Grf/{\rm A}$ ``is" a mapping $Q:Q(1)\to Q(0)$.

The full subcategory $\D$ of $\Grf/{\rm A}$ generated by $\rm D_0$, $\rm D_1$ and $\rm A$ is adequate for $\Grf/{\rm A}$.
Indeed, as explained above, $\Grf/{\rm A}$ can be seen as the category of presheaves on $\D$. 
Furthermore, the inclusions in $\Grf/{\rm A}$ of the full subcategories $\open{\rm A}$ and $\closed{\rm A}$ of left and right mappings 
can be seen as induced by the functors from $\D$ to the arrow category that reduce $\delta_1$ (respectively, $\delta_0$) to the identity.
So $\open{\rm A}$ and $\closed{\rm A}$ are reflective and coreflective.

Now we check that $\Sp{\rm A}$ satisfies the other requirements to be a (strong) bipolar space. $\Grf/{\rm A}$ is a category with tensor, 
the points and the components of $P\tto^{\pi}{\rm A}$ being respectively the arrows, $P(0,1)$, and the components of its ``total" 
graph $P$ (see Proposition~\ref{slice}). Since moreover it is cartesian closed, every part has a complement-negation.
Note that the exponential object $P\imp Q$ can be described by 
\[ (P\imp Q)(0,1) = \Grf/{\rm A}(P,Q) \]  
\[ (P\imp Q)(0) = \Set(P(0),Q(0)) \qq \qq (P\imp Q)(1) = \Set(P(1),Q(1)) \] 
and the domain and codomain mappings are given by $\phi\mapsto \phi(0)$ and $\phi\mapsto \phi(1)$ respectively.

Suppose now that $P\in\open{\rm A}$ and $Q\in\closed{\rm A}$; by the above description of exponential objects and by the definition of 
left and right mappings, it is easy to see that an arrow of $P\imp Q$, that is, a morphism $\phi:P\to Q$, is uniquely determined by 
the corresponding mapping $P(0)\to Q(0)$ which is its domain $\phi(0)$. So $P\imp Q$ is itself a right mapping in $\closed{\rm A}$.
Namely, considering $P$ and $Q$ as mappings $P:P(1)\to P(0)$ and $Q:Q(0)\to Q(1)$, and given $h:P(0)\to Q(0)$, we have
\[ P\imp Q:h\mapsto Q\circ h\circ P:P(1)\to Q(1) \]
and dually $Q\imp P$ is a left mapping in $\open{\rm A}$.
In particular, if $P\in\open{\rm A}$ and $S\in\Set$, the value at $S$ of the functor $\neg P$ is the right mapping
$(\neg P)S$ which acts on mappings $h:P(0)\to S$ by composition: 
\[ (\neg P)S:h\mapsto h\circ P:P(1)\to S \]

To complete the proof that $\Sp{\rm A} = \langle \Grf/{\rm A},\open{\rm A},\closed{\rm A}\rangle$ is a strong bipolar space, 
we have to show that left mappings are coadequate for the right ones, and vice versa; we do this by showing that $\Sp A$ is in fact an atomic space
and so already $\,^{\rm o}(\Sp{\rm A})$ (equivalent to a full subcategory of $\open{\rm A}$) is coadequate for $\closed{\rm A}$ (see Section~\ref{atoms}). 
The parts ${\rm D_0}$ and ${\rm D_1}$ (and only them) deserve to be called (strong) atoms roughly because
\[ \ten({\rm D}_i,P) \cong \hom({\rm D}_i,P) \cong P(i) \qq i=0,1 \]  
To be more precise, it suffices to observe that the unique element in $\ten({\rm D}_i,{\rm D}_i)$ is biuniversal for $\ten$, and so ${\rm D}_i$
is a strong atom. Note that $\up {\rm D_1} = {\rm D_1}$ because it is already a right mapping, while $\up {\rm D_0} = {\rm A}$
because the required universal property of $\delta_0:{\rm D_0}\to {\rm A}$ is assured by the very definition by orthogonality of $\closed{\rm A}$.
Then $(\Sp{\rm A})^{\rm o}$ is the arrow category $\bf 2$ and the associated interpretation takes a right mapping to the corresponding 
covariant presheaf on $\bf 2$. So $(\Sp{\rm A})^{\rm o}$ is adequate for $\closed{\rm A}$ and as a consequence 
(see Corollary~\ref{atomic4}) $\,^{\rm o}(\Sp{\rm A})$ is coadequate for $\closed{\rm A}$ (and symmetrically for left mappings).

Actually, that a (right) functional graph $P$ over $\rm A$ ``is" a mapping $P:P(0)\to P(1)$ in two distinct but equivalent ways 
(Proposition~\ref{atomic3}) is a very intuitive and familiar fact, which may be rephrased as follows: 
\begin{enumerate}
\item
$P$ is the mapping which takes each element of $P(0)$ to the codomain, in $P(1)$, of {\em the} arrow starting from it, 
the definite article being justified by the above definition of (right) functionality.
This may be summarized with the slogan ``follow the arrows".
\item
$P$ is the mapping which takes each element of $P(0)$ to {\em the} element of $P(1)$ in its component, the definite article being justified 
by the fact that components of $P$ are in bijective correspondence with $P(1)$ (see e.g. Proposition~\ref{pcoadj};
indeed, as a category with tensor $\closed{\rm A}$ is particularly simple: $\point P \cong P(0)$ and $\comp P \cong P(1)$). 
This may be summarized with the slogan ``take the components".
\end{enumerate} 
Note that $P$ serves both as the ``graph" and the ``cograph" of the corresponding mapping: 
in the former case, $P$ is seen as the subset of $P(0)\times P(1)$ given by $p_0\wedge p_1:P(0,1)\to P(0)\times P(1)$;
in the latter, it is seen as the quotient of the sum $P(0) + P(1)$ (the set of nodes of $P$) induced by $\eta_P:P\to\disc\comp P$.
In the Appendix, a related definition of ``comapping" will be given.

The proof that $\Sp{\rm A}$ is a strong atomic bipolar space is now complete, and we may apply proposition~\ref{atomic5} to obtain
the following formulas for the reflection and coreflection of $\Grf/{\rm A}$ in right mappings (and dual ones for those in left mappings):
given a graph $P$ on ${\rm A}$, $P\rup$ is the codomain mapping $p_1:P(0,1)\to P(1)$,
while $\up P$ is the mapping $P(0)\to\comp P$ which takes a node to the corresponding component of the total graph.  

As an illustration of the above formulas and of the remarks about exponentiation at the end of Section~\ref{bip}, if $P$ 
is the right mapping such that $P(0)$ and $P(1)$ have two and one element respectively, and if $Q$ is the constant part $\disc S$
on a two element set $S$, the reader may easily compute the left mapping $P\imp Q$ and verify that its coreflection $(P\imp Q)\rup$ 
in right mappings gives the exponential in $\closed{\rm A}$. 

\subsection{Graphs and endomappings}
\label{endomappings}

We now proceed like above, but with the loop graph $\rm L$ instead of the arrow $\rm A$.
Since $\rm L$ is terminal in $\Grf$, $\Grf/\L = \Grf$ and we are going to present the bipolar space
\[ \Sp\L = \langle \Grf,\open\L,\closed\L\rangle \] 
which has left and right endomappings as open and closed parts.
If $\D$ is the full subcategory of $\Grf$ generated by the dot $\rm D$ and the arrow, then $\D$ is adequate for $\Grf$,
which can be seen as the category of presheaves on it: a graph $P\in\Grf$ is determined by its arrows $P(a)$, 
its nodes $P(n)$, and the domain and codomain mappings $p_0,p_1:P(a)\to P(n)$. 

On the other hand, $\D$ is not coadequate. In fact, since $\D\op\cong\D$,
we have a notion of ``dual graph" $P^*$, corresponding to the cointerpretation of a graph $P$ via the inclusion 
of $\D$ in $\Grf$: $P^*(a) = P(n)$, while $P(n) = \comp (P\times {\rm A})$. For example,
the dual of the dot is the arrow, and in general, if $n\geq 1$ then ${\rm C}_n^*\cong {\rm C}_{n+1}$; 
the dual of a right endomapping (see below) is the opposite left endomapping. Note also that any graph and its thin reflection have the same dual graphs,
and that duality preserves sums (in fact any colimit) but not products (as that of the dot with itself).

Again, right and left endomapping are defined by orthogonality with respect to ${\rm D}\tto^{\delta_0}{\rm A}$
and ${\rm D}\tto^{\delta_1}{\rm A}$, and the resulting subcategories $\closed\L$ and $\open\L$ may be seen as induced
by the functors from $\D$ to $\N$ (the monoid of natural numbers) which send $\delta_0$ to $0$ and
$\delta_1$ to $1$, or the other way round; so, by Kan extensions, they are reflective and coreflective in $\Grf$.
Note that, as for $\open{\rm A}$ and $\closed{\rm A}$, $\open\L$ and $\closed\L$ are the ``same" as
categories in themselves (since $\rm A$ and $\rm L$ are self-dual), but they are ``opposite" as subcategories of $\Grf$. 

In $\Grf$, $P\imp Q$ is given by 
\[ (P\imp Q)(a) = \Grf(P\times {\rm A},Q) \] 
\[ (P\imp Q)(n) = \Set(P(n),Q(n)) \] 
and, as in the case of graphs over $\rm A$, given $P\in\open\L$ and $Q\in\closed\L$ we have
that $P\imp Q$ is the right endomapping
\[ P\imp Q:h\mapsto Q\circ h\circ P:P(n)\to Q(n) \]
for any mapping $h:P(n)\to Q(n)$.  
Indeed, recall that the total functor is left adjoint to $- \times A : \Grf\to \Grf/{\rm A}$ 
and so the morphisms $P\times {\rm A}\to Q$ correspond to morphisms
$P\times {\rm A}\to Q\times {\rm A}$ in $\Grf/{\rm A}$, and observe that
the effect of multiplying a right (or left) endomapping by the arrow $\rm A$, is that of displaying it as a
right (or left) mapping in $\Grf/{\rm A}$.

In particular, given a set $S$, $(\neg P)S$ is the right endomapping
\[ (\neg P)S:h\mapsto h\circ P:P(n)\to S \]
for any mapping $h:P(n)\to S$. For example, if $S$ has two elements, $(\neg P)S$ acts on subsets of $P(n)$
like the contravariant power set functor applied to $P$; or if $P = {\rm C}_{\infty}\op$, $(\neg P)S$ is the translation
endomapping on the set of sequences in $S$ (on the other hand, $P\imp^{\open\L}S$ can be obtained again as the 
coreflection of $(\neg P)S$ in $\open\L$ (see Proposition~\ref{clopen}) and is the set of ``bisequences" $Z\to S$). 

As in Section~\ref{mappings}, the dot $\rm D$ is the unique atom of $\Sp\L$, in fact a strong one; but now 
its refection in right endomappings is the infinite chain: $\up {\rm D} = {\rm C}_{\infty}$
(and dually, $\down {\rm D} = {\rm C}_{\infty}\op$).
Then $(\Sp\L)^{\rm o}$ is the monoid $\N$ of natural numbers and the associated interpretation takes a right mapping 
to the corresponding presheaf on $\N$, obtained by its iterations.
So $(\Sp\L)^{\rm o}$ is adequate for $\closed\L$ and as a consequence $\,^{\rm o}(\Sp\L)$ is coadequate 
for $\closed\L$ (and dually for open parts).
That is, a right functional graph $P$ ``is" an endomapping $P(n)\to P(n)$ in two distinct but equivalent ways:
\begin{enumerate}
\item
each node in $P(n)$ ``is" a chain ${\rm C}_{\infty} \to P$, which translated gives its image; 
\item
each node in $P(n)$ ``is" also a ``cochain", that is a component of  ${\rm C}_{\infty}\op \times P$, which translated gives its image. 
\end{enumerate}  

The proof that $\Sp\L$ is a strong atomic bipolar space is so complete, 
and as before we get the following formulas for the reflection and coreflection of $\Grf$ in (right) endomappings:
the coreflection $P\rup$ of a graph $P$ is the endomapping which acts on the set $\Grf ({\rm C}_{\infty},P)$
of the chains of $P$ by translation, while $\up P$ is the endomapping which acts on the set 
$\ten({\rm C}_{\infty}\op,P) = \comp ({\rm C}_{\infty}\op\times P)$
of cochains of $P$ again by translation (each component therein is sent on the component of the translation
of any of its nodes). Of course, the same results could have been obtained, in this case, from the end and the coend formulas for the 
Kan extensions along the above mentioned functor $K:\D\to\N$; indeed, $\ten({\rm C}_{\infty}\op,P)$ is nothing but the tensor
product of the set functors of opposite variance $P$ and $\N(K-,*) = {\rm C}_{\infty}\op$, to which that coend reduces 
in the case of set functors (see Proposition~\ref{coend}). On the other hand, in contexts such as
that of Section~\ref{categories}, the Kan extensions machinery is no longer adequate and the above formulation appears to be the proper one;
in fact we shall partially reverse the usual conceptual hierarchy, deriving the coend formula for the left Kan extensions of set functors
from general properties of the ``bipolar doctrine" $\Sp X , X\in\Cat$ in Section~\ref{kan}.

We end this section with some examples of reflection and coreflection of graphs on endomappings which are easily
``computed" with the above formulas (for the reflections, it may be useful to sketch on a piece of paper a portion
of the product ${\rm C}_{\infty}^{\rm op}\times P$ to understand how $\up P$ works):
\begin{enumerate}
\item \label{ex1}
Among the chains, the only (right) endomappings are ${\rm C}_0$ (the empty graph) and ${\rm C}_{\infty}$:
the ${\rm C}_n$ are not endomappings because they have a node which is not the domain of any arrow.
If $1\leq n<\infty$ then ${\rm C}_n\rup = {\rm C}_0$ and $\up{\rm C}_n = {\rm C}_{\infty}$.
\item \label{ex2}
To examine the other case of non functionality, that is the presence of nodes which are domain of more then one arrow,
let ${\rm S}_n$ be the graph obtained by the discrete $\disc S$ on an $n$-elements set $S$, by adding
a further node $s$ and $n$ arrows, all starting from $s$ and ending on each of the $n$ loops; note that only ${\rm S}_1$
is an endomapping. Then, if $n\geq 1$, ${\rm S}_n\rup$ is the copower $n{\rm S}_1$, while $\up{\rm S}_n$ is ${\rm S}_1$.
\item
Let $P$ be a graph with just one node; then $P\rup$ is the endomapping which acts on 
the set of sequences in its set of arrows by translation;
similarly, if $P$ is codiscrete $P\rup$ acts on the set of sequences in its set of nodes.   
In both cases, $\up P$ is the terminal graph $\L$.
\item
If $P$ is obtained by adding a loop to the first node of ${\rm C}_{\infty}$, then $P\rup$ is the double infinite
chain $\dch$; $\up P$ is again the terminal graph. 
\end{enumerate}
As it is clear from examples~\ref{ex1} and~\ref{ex2} above, $(-)\rup$ tends to act ``on the domains", deleting those nodes which have no arrows 
out of them and duplicating those which have more than one. On the other hand, $\up(-)$ tends to act ``on the codomains", creating those nodes 
which are needed and collapsing those which appear as ends of arrows out of the same node.

\section{Subsets and sieves}
\label{preorders}

In this section we present the two-valued bipolar space $\Sp X = \langle {\cal P}(X),\open X,\closed X\rangle$ associated to a poset $X$,
which is the two-valued version of the main exemple to be studied in the next section. 
The connection with Alexandrov spaces is considered, and it is shown how the formulas for reflections and coreflections of parts are also
suggestive for general topological spaces.

Recall that, given a poset $X$, a sieve or lower set of $X$ is a $\Two$-valued presheaf $A:X\op\to\Two$, 
or equivalently a subset $A$ of $X$ such that $x\in A$ and $y\leq x$ implies $y\in A$;
cosieves (or upper sets) are defined dually. 
Sieves and cosieves form the subposets $\open X$ and $\closed X$ of $({\cal P}(S),\subseteq)$, the set of parts of the underlying set of $X$;
being presheaves categories, $\open X$ and $\closed X$ are complete lattices.

To show that $\Sp X = \langle {\cal P}(X),\open X,\closed X\rangle$ is indeed a strong bipolar space it suffices to observe that 
\begin{enumerate}
\item
${\cal P}(X)$ has tensor, with $P\meets Q$ iff $P$ meets $Q$, that is iff they have a non empty intersection.
\item
Since a part of $X$ is a sieve iff its complement in ${\cal P}(X)$ is a cosieve, each open part has a closed complement,
and vice versa.
\item
The complement is classical, since $\P(\Sp X) = {\cal P}(X)$ is in fact a Boolean algebra.
\item
Since the inclusion of the discrete poset with the same elements (objects) of $X$ in $X$ itself induces the inclusion of $\open X$
in ${\cal P}(X)$, by Kan extensions the former is reflective and coreflective in the latter, and similarly for $\closed X$.
\item
Since, given $A\in\open X$ and $D\in\closed X$, the relative complement $A\imp D$ in the Boolean algebra ${\cal P}(X)$
is given by $D\cup(\neg A)$, the bipolar space $\Sp X$ is strong.
\end{enumerate}
Furthermore, the space $\Sp X$ is atomic; indeed, the parts $x$ with a single element are strong atoms: 
\[ x\meets P\iff x\subseteq P\iff x\in P \]  for any part $P$, and their reflections $\down x$ and $\up x$ in $\open X$ and in $\closed X$ 
are the principal (co)sieves, or representable presheaves:
\[ y\in\down x\iff y\leq x \qq\qq  y\in\up x\iff x\leq y \]
Then the strong atoms are the objects of subposets of $\,^{\rm o}(\Sp X)$ and $(\Sp X)^{\rm o}$ isomorphic to $X$ and $X\op$ respectively, and are 
manifestly adequate for $\open X$ and $\closed X$. Note that, if $X$ is not partially ordered (i.e. skeletal), there are non-strong atoms, namely 
the parts formed by pairwise isomorphic elements-objects. 

Observe also that any two-valued bipolar space that has ${\cal P}(X)$ as parts, is of the form $\Sp X$ for a preorder on $X$.
Indeed, the strong atoms are again the elements of $X$ as parts in ${\cal P}(X)$, and the corresponding subposets of $\oX$ and $\Xo$ 
as above define the desired preorder: 
\begin{enumerate}
\item
any open part $A\in\open X$ is a sieve, namely that given by its interpretation (and similarly for closed parts);
\item
any sieve $A\subseteq X$ is the union of the open reflections $\down x$ of its atoms $x\in A$, and so is itself open.
\end{enumerate}
The open parts of $\Sp X$ form a topology on the underlying set of $X$, namely the so-called Alexandrov topology associated to the poset $X$, 
which has closed parts as closed sets.
Actually, the hall-mark of such topological spaces is the reflectivity of its open sets (or dually, the coreflectivity of the closed ones): 
every part $P\in{\cal P}(X)$ has not only an ``interior" and a ``closure" but also an open reflection $\down P$ and a closed coreflection $P\rup$. 
Let us summarize the above considerations:
\begin{prop}
The bipolar spaces $\Sp X$ associated to a poset on $X$ are coextensive with bipolar spaces with ${\cal P}(X)$ as parts and with Alexandrov
topological spaces on $X$.
\end{prop}
\epf
\noindent The definition of continuous maps between bipolar spaces, along with the appropriate extension of the 
above proposition will be carried out in Section~\ref{doctrines}. 

Formulas~(\ref{refp}) and~(\ref{corp}) of Section~\ref{two}, which give the reflections and the coreflections of a part $P$ of an
atomic two-valued bipolar space, then apply to Alexandrov spaces, but not to general topological spaces. Yet, the first of the~(\ref{refp}) 
and of the~(\ref{corp}) still have a meaning if, as a surrogate of the open reflection $\down x$, the neighborhood filter of a point $x\in X$ is considered.
Indeed, thinking of $\down x$ as the set of ultrafilters converging to $x$, 
\[ x\in\up P \iff \down x\meets P \]
can be rephrased by saying that $x$ is in the closure of $P$ iff its neighborhood filter meets $P$, that is, if there is an ultrafilter in $P$ 
converging to $x$, and similarly 
\[ x\in P\rdown \iff \down x\subseteq P \]
can be rephrased by saying that $x$ is in the interior of $P$ iff its neighborhood filter is included in $P$, that is, if any ultrafilter 
converging to $x$ is in $P$.

\section{Categories over a category}
\label{categories}

This and the next sections are devoted to the situation that motivated the notion of bipolar space itself: 
the simultaneous inclusion of discrete fibrations and discrete opfibrations on a category 
in the category of categories over it.
Indeed, Sections~\ref{graphs} and~\ref {preorders} displayed bipolar spaces $X = \langle \P(X),\open X,\closed X\rangle$ 
in which $\P(X)$ is a presheaf category $[\A,\Set]$ and the inclusions of $\open X$ and $\closed X$ in it are induced by functors from $\A$ 
(in fact, connected functors, since the inclusions are full, see~\cite{elv02}), and so could have been considered in the frame of Kan extensions;
presently we will study a case in which, on the contrary, $\P(X)$ is not even cartesian closed (although $\open X$ and $\closed X$ 
are still presheaf categories) and for which the new conceptual tools are more appropriate.

After briefly recalling some well-known basic facts about the slice categories $\Cat/X$, particularly concerning exponentiability therein,
we present the bipolar space 
\[ \Sp X = \langle \Cat/X,\open X,\closed X\rangle \] 
already mentioned in the Introduction, and show how some important categorical notions turn out to be enclosed in it; 
indeed, $\Sp X$ usefully extends a fragment of the bicategory of modules (or bimodules, or profunctors) associated to the category $X$.  
Next, we turn to atoms in $\Sp X$, considering first those associated to the objects of $X$ and then the more general ones
associated to the idempotent arrows in $X$, obtaining in particular the formulas for the reflection and the coreflection in open and closed parts 
(that is, discrete fibrations and discrete opfibrations) of a part of $\Sp X$ (that is, a category over $X$).
 
\subsection{Categories over a category and exponentiability}
\label{catexp}

$\Cat$ is a category with components. Indeed, the points $\point X$ of a category $X$ are its objects,  
its components $\comp X$ are the components of the underlying graph, and the discrete categories are those whose only arrows are the identities.
So, being cartesian closed, $\Cat$ is a category with tensor, and the properties of Proposition~\ref{closed} hold for categories.
Then, by Proposition~\ref{slice}, every slice $\Cat/X$ has tensor, components being given by those of the corresponding total categories.

The slice over the arrow category $\2$ is specially meaningful. The objects $P\tto^{\pi}\2$ of $\Cat/\2$ can be identified 
with modules: two categories $P_0$ and $P_1$ which act on a ``matrix" or sets $P(x,y)$ ($x\in P_0$, $y\in P_1$), on the left and on the right respectively, 
in such a way that the associative law 
\[ (ua)v = u(av) \] 
holds ($u:x'\to x$ in $P_0$, $v:y\to y'$ in $P_1$ and $a\in P(x,y)$).
The morphisms $\phi : P\to Q$ in $\Cat/\2$ are functors $\phi_0:P_0\to Q_0$ and $\phi_1:P_1\to Q_1$, and a family of mappings
$\phi(x,y):P(x,y)\to Q(x,y)$ such that: 
\[ \phi(x',y')(uav) = (\phi_0 u)(\phi(x,y)a)(\phi_1 v) \]

Among the objects of $\Cat/X$, there are the subterminal ones corresponding to the objects and to the arrows of $X$; e.g., if $x\in X$,
we denote by $x$ also the object of $\Cat/X$ given by the corresponding functor $x:\1\to X$ (and similarly for arrows $f:\2\to X$).
Note that if $P\tto^{\pi}X$ is an object over $X$, multiplying it by $x$ in $\Cat/X$ gives (as total category) 
the fibre category $Px$ over $x$, while multiplying it by $f:x\to y$ gives the module $Pf$ over $f$, with $(Pf)_0=Px$ and $(Pf)_1=Py$.
So, an object $P\tto^{\pi}X$ of $\Cat/X$ is in particular a family of categories, the fibres $Px$ ($x\in X$), and a 
family of modules, the fibres $Pf$ ($f$ in $X$); likewise, a morphism $\phi : P\to Q$ in $\Cat/X$ is determined by the corresponding
families of functors between the fibres $\phi(x):Px\to Qx$ and of module morphisms $\phi(f):Pf\to Qf$ ($f:x\to y$ in $X$), 
with $\phi(f)_0 = \phi(x)$ and $\phi(f)_1 = \phi(y)$.

While $\Cat$ is cartesian closed, it is not locally cartesian closed, that is not all its slices ${\Cat}/X$ are such.
If $f:x\to y$ and $g:y\to z$ are consecutive arrows in $X$, seen as objects of $\Cat/X$, then their composite $h=g\circ f$ (as a functor $\3\to X$) 
is also their pushout with respect to $y$ in $\Cat/X$.
Then, if $P$ is exponentiable, $Ph=P\times h$ should be the pushout of $Pf=P\times f$ and $Pg=P\times g$, and the latter
is easily seen to be given by their composite as consecutive modules $Pg\otimes Pf$ (with the usual coend formula, see Section~\ref{disfib});
thus a necessary condition for the exponentiability is that 
\[ Pg\otimes Pf \cong P(g\circ f) \] 
(that is, the morphisms $Pg\otimes Pf \to P(g\circ f)$ given by composition in $P$ are actually isomorphisms)
which, explicited, becomes the well-known factorization lifting property (see e.g.~\cite{joh99}, \cite{bun00}, \cite{str01} and references therein). 

To verify that the condition is also sufficient, observe that, analyzing the expected $P\imp Q$ using figures whose shapes 
are the $x$ or the $f$ as above, we deduce that $(P\imp Q)x = \Cat(Px,Qx)$ and the arrows over $f$ in $(P\imp Q)f$ 
are the module morphisms $\phi:Pf\to Qf$, the domain and codomain being given by the corresponding $\phi_0$ and $\phi_1$. 
If $\phi$ and $\psi$ are consecutive arrows of the expected $P\imp Q$,
over $f$ and $g$ respectively, then the above condition allows a meaningful definition of composition 
\[ \psi\circ\phi : P(g\circ f) \to Q(g\circ f) \]
by the very property of the pushout $P(g\circ f) = Pg\otimes Pf$.

For instance, if $X=\3$, the inclusion $\iota :\2\to\3$ which preserves the initial and the terminal objects is not 
exponentiable in $\Cat/X$, since, as can be easily checked, $\iota\imp \disc S$ does not exist, if the set $S$ has more than one element.

\subsection{Discrete fibrations and negation}
\label{disfib}

Recall that a category $P\tto^{\pi}X$ on $X$ is said to be a discrete opfibration (abbreviated, dof) if, given any arrow $f:x\to y$ in $X$ 
considered as an object $\2\to X$ of $\Cat/X$, it is interpreted as a bijection by the corresponding domain arrow $x\to f$ in $\Cat/X$.
Discrete fibrations (df) on $X$ are defined dually. 
We denote by $\open X$ and $\closed X$ the full replete subcategories of $\Cat/X$ generated by the df's and, respectively, the dof's on $X$.
By a well-known argumentation, these are equivalent to the presheaf categories $[X\op,\Set]$ and $[X,\Set]$ respectively.
In a direction, one associates to a presheaf on $X$ its category of elements; in the other, given a df $A$ over $X$, every fibre $Ax$
is a discrete category which gives the value of the corresponding presheaf at $x$, and for any arrow $f:x\to y$ in $X$ and $b\in Ay$, $(Af)b$ (or simply $fb$) 
will denote the domain of the unique lifting of $f$ with codomain $b$. 
When appropriate, we will not distinguish notationally a presheaf $A$ and the corresponding df in $\Cat/X$, 
e.g. $Ax$ will denote both the value of the presheaf $A$ at $x$ and the fibre on $x$ of the corresponding df.

Satisfying the factorization lifting condition of Section~\ref{catexp}, the df's and the dof's are exponentiable, 
and so are the products of a df and a dof in $\Cat/X$
(in fact all of them belong to the special class of UFL functors, see~\cite{bun00}).
If $A$ is a df and $P\tto^{\pi}X$ is an object of $\Cat/X$, then, by the above description of exponentials, $(A\imp P)x$
is the set of mappings 
\[ h:Ax\to\point Px \] 
from the set $Ax$ to the set of the objects of $Px$; 
if furthermore $D$ is a dof, then the modules $Af$ and $Df$ ($f$ in $X$) are essentially
left and right mappings as in Section~\ref{mappings}. So, like in the case graphs over the arrow, we have 
\begin{prop}   
If $A$ is a df and $D$ is a dof on $X$, then $A\imp D$ is a dof which, as a covariant presheaf, acts as follows:
given a mapping $h:Ax\to Dx$ in $(A\imp D)x$ and an arrow $f:x\to y$ in $X$, 
\[ (A\imp D)f : h\mapsto Df\circ h\circ Af : Ay\to Dy \]
and dually $D\imp A$ is a df in $\open{\rm X}$.
\end{prop}
\epf
\noindent In particular, since any discrete object (that is, a constant presheaf) $\disc S$ of $\Cat/X$ is both a df and a dof,
\begin{corol} \label{negcor}
Any $A$ in $\open X$ has a negation functor $\neg A:\Set\to \Cat/X$, right adjoint to $\ten(A,-)$, which takes values in closed parts. 
The value of $\neg A$ at $S\in\Set$ is the exponential $A\imp\disc S$ that, considered as a (covariant) presheaf, is $\Set(A-,S)$. 
Dually, any dof $D$, has  a negation $\neg D$ with open values.
\end{corol}
\epf
\noindent Explicitly, $(\neg A)S:X\to\Set$ acts as follows:
\begin{eqnarray}   
x\mapsto\Set(Ax,S)   \nonumber               \\
f:h\mapsto h\circ Af : Ay\to S  \label{negc} 
\end{eqnarray}
for any arrow $f:x\to y$ in $X$ and any mapping $h:Ax\to S$.  
Note that an object in $\Cat/X$ may not have negation, as the above example with $X=\3$ shows. 

Then, observing that the clopen parts of $\Sp X$ are the discrete bifibrations, which correspond to presheaves that act by bijections, 
Corollary~\ref{clopencor} applies, giving 
\begin{corol}
Discrete bifibrations form a cartesian closed full subcategory $\clopen X$ of $\Cat/X$, and the inclusion functor preserves exponentials.
In particular, if $X$ is a grupoid, the open (or closed) parts of the corresponding Boolean space $\Sp X$ form a cartesian closed full 
subcategory $\open X = \closed X$ of $\Cat/X$, wherein exponentials are computed as in $\Cat/X$.
\end{corol}
Thus we have a conceptual explanation of the fact that ``group theory is simpler than category theory" (\cite{law70}), in that the Frobenius law
hold in the grupoid doctrine (see Section~\ref{doctrines}).

\subsection{The tensor functor and the tensor product}
\label{modules}

Now we show that $\Cat/X$ encompasses a small portion of the bicategory $\Mod$, whose $0$-cells and $1$-cells are categories and modules between them,
while the $2$-cells can be identified with those module morphisms of Section~\ref{catexp}  
which leave the domain and the codomain unchanged (that is $\phi_0$ and $\phi_1$ are identity functors). 

Recall that if $R:X\to Y$ and $S:Y\to Z$ are consecutive modules in $\Mod$, their composition $S\otimes R : X\to Z$ 
is given by the coend 
\[ (S\otimes R)(x,z) = \int^y R(x,y)\times S(y,z) \]
while if $R:X\to Y$, $S:X\to Z$ and $T:Z\to Y$ then the biclosed structure of $\Mod$ gives the modules $R\triangleright S:Y\to Z$ 
and $T\triangleleft R:Z\to X$ as the ends 
\begin{eqnarray} \label{bic}
(R\triangleright S)(y,z) =\int_x \Set(R(x,y),S(x,z)) \\ 
(T\triangleleft R)(z,x) =\int_y \Set(R(x,y),T(z,y))
\end{eqnarray}
It is well-known that the structure of $\Mod$ (and in particular the fact that it is biclosed) 
encodes important aspects of category theory (such as generalized Kan extensions, see e.g.~\cite{law73}). 
We concentrate on a fragment of it, namely that generated by a category $X$ and by $\1$, the terminal category.
In that fragment we find in particular the category $\Mod(\1,\1) = \Set$ of ``truth values", 
the category $\Mod(X,X) = [X\op\times X,\Set]$ of endomodules on $X$, and the presheaf categories 
\[ \Mod(\1,X) = [X\op,\Set] \qq\qq \Mod(X,\1) = [X,\Set] \]
The composition bifunctor $\otimes$ in $\Mod$ then gives in particular 
\[ \otimes : [X,\Set] \times [X\op,\Set]\to\Set \] which is the well known
tensor product between presheaves of opposite variance (see also~\cite{mam91}); 
if $A\in {\Mod}(\1,X)$ and $D\in \Mod(X,\1)$, the set $D\otimes A$ is given by the coend of the functor $A\otimes D\in\Mod(X,X)$:
\[ D\otimes A = \int^x A\otimes D = \int^x Ax\times Dx \]
Furthermore, its closed structure gives adjunctions 
\[ D\otimes -\adj -\triangleleft D :\Set\to [X\op,\Set] \qq\qq  -\otimes A\adj A\triangleright - :\Set\to [X,\Set] \] 
with parameters $D\in [X,\Set]$ and $A\in [X\op,\Set]$.

\begin{prop}  \label{coend}
The functor $\ten: \Cat/X\times \Cat/X \to \Set$ extends the functor $\otimes : [X\op,\Set] \times [X,\Set]\to\Set$;
that is, given a df $A$ and a dof $D$ on $X$, 
\[ \ten(A,D) = D\otimes A \] 
the $A$ and $D$ on the right side standing for the corresponding set functors.
Furthermore, the negation $\neg A$ of a df $A$ is the same as $A\triangleright -$ and dually, $\neg D = -\triangleleft D\,$. 
\end{prop}
\pf Indeed, the colimit in $\Set$ corresponding to the coend $\int^x Ax\times Dx$ of the bifunctor $A\otimes D:X\op\times X\to\Set$ 
(see~\cite{mam91}) is easily seen to reduce to the components of the category whose objects are the pairs $(a,b)\in Ax\times Dx$
and whose arrows from $(a,b)\in Ax\times Dx$ to $(a',b')\in Ay\times Dy$ are the arrows $f:x\to y$ in $X$ such that $Af:a'\mapsto a$ and $Df:b\mapsto b'$;
then the first part of the proposition follows from the fact that this category is the product in $\Cat/X$ of the df and the dof 
corresponding to $A$ and $D$.
The second part follows from the fact that, having closed values, $\neg A$ serves as a right adjoint of $A\otimes - :\closed X\to\Set$
as well; alternatively, observe that equation~(\ref{negc}) coincides with the formula for $S\triangleleft A$ of equation~(\ref{bic}), 
since in this case the universal quantification (the end) is on the terminal category. 
\epf
Of course, if $A,B\in\open X$ and $D,E\in\closed X$, we also have
\[ \hom(A,B) = B\triangleleft A \qq\qq \hom(D,E) = D\triangleright E \] 
since e.g. the end $D\triangleright E = \int_x \Set(Dx,Ex)$ is the set of natural transformations $D\to E$.

\subsection{Objects as atoms and their reflections} 
\label{atomcat}

To show that $\Sp X$ is indeed a (strong) bipolar space, we have to check that the negation
is classical, and that $\open X$ and $\closed X$ are reflective and coreflective in $\Cat/X$.
For the former fact, we show presently that, as in the other examples, $\Sp X$ is in fact an atomic space;
for the latter we presently draw the form that the expected reflections and coreflections should have (Proposition~\ref{main} 
and Corollary~\ref{maincor}), while in Section~\ref{reflections} we shall prove that these formulas actually work properly.

As mentioned in the introduction, the objects $x:\1\to X$ of $X$, as categories over $X$, are atoms of $\Sp X$; indeed, for a part $P\in\Cat/X$ 
\[ \hom(x,P) \cong \disc(Px) \qq\qq \ten(x,D) \cong \comp(Px) \]
so that $x$ ``is included" in the part $P$ in as many ways as are the objects of the fibre $Px$, while $x$ ``meets" $P$ in as many ways as are
its components. In particular, for discrete fibrations $A\in\open X$ and discrete opfibrations $D\in\closed X$, 
\begin{equation} \label{atomc}
\hom(x,A) \cong \ten(A,x) \cong Ax   \qq \qq  \hom(x,D) \cong \ten(x,D) \cong Dx
\end{equation}
and these bijections are induced by the unique element $u\in\ten(x,x)$ which is therefore biuniversal, as required.

We will see in the next section that there may be also other atoms, but in the meantime we concentrate on those 
corresponding to the objects of $X$. The first important fact is that they have well-known reflections, namely the representable functors: 
\[ \down x = X/x = X(-,x) \qq\qq \up x = x/X = X(x,-) \]
Indeed, this is the content of the Yoneda lemma, the units $\lambda_x : x \to X/x$ and $\rho_x : x \to x/X$ of the reflections sending 
the only object of $x\in\Cat/X$ to the object $x\tto^{id} x$ of $X/x$ and of $x/X$ respectively.
Those units then induce bijections
\[ \hom(\down x,A) \cong \hom(x,A) \qq\qq \hom(\up x,D) \cong \hom(x,D) \]
and so by the~(\ref{atomc}) also bijections
\[ \hom(\down x,A) \cong Ax \qq\qq \hom(\up x,D) \cong Dx \]
and (addendum to the Yoneda lemma) these bijections are natural in $x\in X$. This can be rephrased by saying that $A$ is interpreted
as itself by the Yoneda functor $Y:X\to\Cat/X$; or better, the discrete fibration $A$ is interpreted as the corresponding presheaf on $X$.
Symmetrically, any discrete opfibration $D$ is interpreted as the corresponding presheaf on $X\op$ by $Y':X\op\to\Cat/X$.
So, the atoms associated to the objects $x\in X$ are adequate for open and closed parts, and $\Sp X$ is atomic.
From Proposition~\ref{atomic3}, we get the ``co-Yoneda lemma":
\begin{equation} \label{coyo}
\ten(A,\up x) \cong Ax  \qq\qq \ten(\down x,D) \cong Dx
\end{equation}
for any $A\in\open X$ and $D\in \closed X$, with the corresponding addendum that these bijections are natural in $x\in X$;
of course, this is a particular case of the fact that if $R:X\to Y$ and $S:Y\to Z$ are consecutive modules in $\Mod$ and $R$ is represented
by a functor $F:X\to Y$, then $(S\otimes R)(x,z) \cong S(Fx,z)$.

Thus, not only $\Sp X$ is an atomic space, but the objects-atoms $x\in X$ generate full subcategories of $\,^{\rm o}(\Sp X)$ 
and $(\Sp X)^{\rm o}$ which are already adequate for closed and open parts, and which are isomorphic to $X$ and $X\op$ respectively. 
So Proposition~\ref{atomic5} is still valid when restricted to objects-atoms, that is when the Yoneda embeddings 
\[ \down(-):X\to\Cat/X \qq\qq \up(-):X\op\to\Cat/X \]
are used instead of $\lp$ and $\rp$, giving
\begin{prop} \label{main}
The reflections of an object $P\tto^{\pi}X$ of $\Cat/X$ in discrete fibrations and discrete opfibrations are given by:
\[ \down P = \ten(P,\up -) \qq \qq \up P = \ten(\down -,P) \]
and the coreflections by
\[ P\rdown = \hom(\down -,P) \qq \qq P\rup = \hom(\up -,P) \]
That is, they are obtained as the cointerpretation and interpretation of $P$ via the Yoneda embeddings.
\end{prop} 
\epf
\noindent Note that we are following the common ``abuse" of denoting the object $P\tto^{\pi}X$ over $X$ simply with $P$.

It is immediate to see that multiplying $\up x$ or $\down x$ with $P$ in $\Cat/X$ one obtains the ``comma" or ``map" categories
\[ \up x\times P = x/P  \qq\qq  \down x\times P = P/x \] 
where $P/x = \pi/x$ has as object the pairs $(a,f)$ with $a\in P$ and $f:\pi a\to x$ in $X$, and as arrows $(a,f)\to (b,g)$ the 
arrows $u:a\to b$ in $P$ such that $f = g\circ\pi u$; $x/P$ is defined dually.
So we find again the result (see~\cite{par73} and~\cite{law73}):
\begin{corol} \label{maincor}
The reflections of an object $P\tto^{\pi}X$ of $\Cat/X$ in discrete fibrations and discrete opfibrations are given by
the components of the corresponding map categories
\[ \down P = \comp(-/P) \qq \qq \up P = \comp(P/-) \]
\end{corol} 
\epf
Alternatively, to check that discrete fibrations and discrete opfibrations are reciprocally coadequate in $\Cat/X$,
one might prove directly the contraposition law, as we now briefly sketch.
A morphism $\alpha:A\to B$ in $\open X$ is a natural transformation $\alpha(x):Ax\to Bx$ and, by Corollary~\ref{negcor}, its 
negation $\neg\alpha:\neg B\to\neg A$ is the natural transformation between the corresponding negation functors $\Set\to\closed X$ 
defined by the family of mappings
\[ \neg\alpha(x,S) = \Set(\alpha(x),S):\Set(Bx,S)\to \Set(Ax,S) \] 
On the other hand, any natural transformation $\Theta:\neg B\to\neg A$ is expressed by a family of 
morphisms $\Theta(S):(\neg B)S\to(\neg A)S$, and so ultimately by a family of mappings $\Theta(x,S):\Set(Bx,S)\to\Set(Ax,S)$, 
natural in $x\in X$ and $S\in\Set$. For any $x\in X$, by the naturality in $S$ and by Yoneda we have 
a mapping $\alpha_\Theta(x):Ax\to Bx$. The family $\alpha_\Theta(x)$ is also natural in $x\in X$ and so defines a 
morphism $\alpha_\Theta:A\to B$ whose negation is $\Theta$. By Yoneda again, such a morphism is uniquely determined.  

\subsection{Idempotents as atoms and the Cauchy completion} 
\label{atomcat2}

We now show that to any idempotent arrow $e:x\to x$ in $X$ there corresponds an atom in the bipolar space $\Sp X = \langle \Cat/X,\open X ,\closed X\rangle$,
whose open and closed reflections are retracts of the functors represented by $x$. The results on object-like atoms of Section~\ref{atomcat} can 
be obtained specializing those of the present section to the case in which the idempotent mapping $e:x\to x$ is the identity $\id_x$.

To begin with, let us recall the basic notion of Isbell conjugation (see~\cite{law86}), and fix some terminology.
Given a category $X$,there are adjoint functors
\[ (-)^*\adj (-)^\# : [X,\Set]\op\to [X\op,\Set] \]
where the conjugate of $A:X\op\to\Set$ is the functor $A^*:X\to\Set$ defined on objects $x\in X$ by
\[ A^* x = [X\op,\Set](A,\down x) \]
and dually for the conjugate $D^\#$ of $D:X\to\Set$:
\[ D^\# x = [X,\Set](D,\up x) \]
More precisely, $A^*$ is the interpretation of $A\in [X\op,\Set]$ via the opposite of the Yoneda embedding $Y\op : X\op \to [X\op\to\Set]\op$, 
and similarly for $D^\#$.

The pairs $\langle A,D \rangle$, with $A:X\op\to\Set$ and $D:X\to\Set$, such that $A \cong D^\#$ and $D\cong A^*$,
are called {\bf Dedekind cuts}, since for a poset $X$ they give its Dedekind-MacNeille completion.
Coming back to the study of the bipolar space $\Sp X$, we have the following
\begin{prop}  \label{ded}
Any atomic pair $\langle \down e,\up e\rangle \in \open X \times\closed X$ is a Dedekind cut. 
\end{prop}
\pf Indeed, by the very definition of atomic pair (see Section~\ref{atomic}), 
\[ \ten(\down e, \up x) \cong \hom(\up e,\up x) \qq \qq \ten(\down x,\up e) \cong \hom(\down e,\down x) \]
both the bijections being mediated by a biuneversal element $u\in\ten(\down e, \up e)$.
The result now follows from the co-Yoneda lemma (see equations~(\ref{coyo}) of Section~\ref{atomcat}).
\epf
Note that we are using again the notations of the previous subsections: $\up x$ and $\down x$ are the representable functors, 
and the $\hom$ is taken in $\Cat/X$, so that if $A$ and $B$ are both df's or both dof's, then $\hom(A,B)$ can be thought of as the set 
of natural transformations between the corresponding functors.

Now we turn to the typical non-object atoms which $\Sp X$ may have.
Let $\bf e$ be the category with an object and just one non-identity arrow which is idempotent, so that $\Cat({\bf e},X)$ represents the idempotent 
arrows of $X$. Any idempotent arrow $e:x\to x$ in $X$, as a category ${\bf e}\tto^e X$ over $X$, is an atom of $\Sp X$:
\begin{prop}  \label{idat}
Given an idempotent arrow $e:x\to x$ in $X$, the unique element $u\in\ten(e,e)$ induces bijections
\[ \hom(e,A) \iso \ten(A,e) \qq\qq \hom(e,D) \iso \ten(e,D) \]
for any $A\in\open X$ and $D\in\closed X$. Furthermore, $\hom(e,D) \cong \ten(e,D)$ is isomorphic to the set $D'x\subseteq Dx$ of elements fixed by 
the mapping $De:Dx\to Dx$, and similarly for open parts.
\end{prop}
\pf
We begin with some remarks:
\begin{enumerate}
\item
In general, an object $a$ over $x$ of (the total category of) a part $P\in\Cat/X$ 
corresponds to a morphism $a:x\to P$ in $\Cat/X$ and so we can define $[a]:=\comp a :\comp x \to \comp P$; since $\comp x = 1$
we so obtain an element $[a]$ of  $\comp P$: ``the component of $a$". In the other direction, any component $u\in\comp P$ is the component 
of (at least) one object. Moreover, if $\alpha:P\to Q$ in $\Cat/X$, then $\comp\alpha:[a]\mapsto [\alpha a]$. 
So, denoting by $\tot:\Cat/X\to\Cat$ the total functor, we have a natural transformation
\[ [-]:\point^\Cat\circ\tot\to\comp:\Cat/X\to\Set \]
whose components are onto mappings.
\item
Since any arrow $\alpha:e\to D$ in $\Cat/X$ corresponds to an object $a\in Dx$ over $x$ and a loop at $a$ over $e$, $\hom(e,D)$ 
is the set of fixed points of the mapping $De:Dx\to Dx$. On the other hand, $\ten(e,D)$ is the set of the components of the subcategory $D_e$ 
of $D$, formed by the fibre $Dx$ and by the arrows over $e$; in fact, apart from the identities, $D_e$ is the graph of the endomapping $De$
as in Section~\ref{endomappings}. 
In general, two elements belong to the same component of an endomapping if they are sent to the same element
by suitable iterations of it; for an idempotent mapping, it is enough to check if they are sent to the same element by the mapping itself:
\[ [a] = [b] \iff ea = eb \]
\end{enumerate}
Denoting with $x$ also the unique object of $\down e\,\times\!\up e$, which is isomorphic to ${\bf e}\tto^e X$ in $\Cat/X$, let us check the 
right universality of $u=[x]$ (i.e. for closed parts; left universality follows by symmetry): $\beta\mapsto\ten(e,\beta)u$ is the mapping 
\[ \hom(e,D) \to \ten(e,D):a\mapsto [a] \]
which takes each object-element $a\in Dx$ fixed by $De$, to its component in $\comp D_e$,
and it is indeed a bijection, since the mapping 
\[ \ten(e,D)\to \hom(e,D):[a]\mapsto (De)a \] 
which takes a component in $\comp D_e$ to the fixed point $(De)a$, where $a$ is any of its objects-elements, is well defined 
(by the above remark) and is readily verified to be its inverse.
\epf
As a consequence, $\langle\down e,\up e\rangle$ is an atomic pair (see Proposition~\ref{atomic1}).
It is worth to see what $\up e$ and $\down e$ are concretely. Observe first that, by applying Corollary~\ref{maincor} to any 
part $P\tto^\pi X$ in $\Cat/X$, we have the formula $(\up P) x = \comp(P/x)$ in which
the (total) category $P/x$ is a df, when considered as a category over $P$ rather than over $X$; indeed, as a presheaf on $P$ it is 
given by $X(\pi-,x)$, so that $\comp(P/x)$ may be seen as the components of this presheaf (that is, its colimit). 
In the case of ${\bf e}\tto^e X$, associated to the idempotent $e:x\to x$ in $X$, we then have that $(\up e) y = \comp(e/y)$ is given by 
the components of the corresponding presheaf on $\bf e$, that is of the idempotent endomapping $X(x,y)\to X(x,y)$ obtained by multiplying by $e$;
but then, as just observed in Proposition~\ref{idat}, the components correspond to fixed points and we have proved the first part of the following 
\begin{prop} \label{k} 
The right reflection $\up e$ of the idempotent atom associated to $e:x\to x$ is the covariant presheaf defined as follows:
$(\up e)y$ is the set of the arrows $f:x\to y$ in $X$ such that $f\circ e = f$, and if $g:y\to z$ then $(\up e)g:f\mapsto g\circ f$.
Furthermore $\up e$ is a retract of $\up x$. Dual statements hold for $\down e$.
\end{prop}
\pf As for the second part, one readily checks that $-\circ e:\up x\to\up e$ is a retraction of the inclusion $\up e\to\up x$.
Then $\up e$ is a retract in $[X,\Set]$, and the associated idempotent $\up x\to\up x$ corresponds to $e$ via the Yoneda embedding.
\epf
\noindent Note that if the idempotent $e:x\to x$ splits already in $X$ as $e=i\circ r$, with $r\circ i=\id_{x'}$, then $\up e$ is isomorphic to the
representable $\up x'$, since the natural transformations $-\circ i:\up e\to\up x'$ and $-\circ r:\up x'\to\up e$
are each other's inverse. So the reflections of the atoms associated to split idempotents are in fact already reflections of objects-atoms.

Let us call an arrow $e:x\to x$ in $X$ ``eventually idempotent" if there is a natural number $n_0$ such such that $e^n=e^{n_0}$ for any $n\geq n_0$.
As for idempotent arrows, any eventually idempotent arrow in $X$ can be seen as a category over $X$, and
argumentations similar to those of Proposition~\ref{idat} easily show that it is an atom of $\Sp X$. Thus $\Sp X$ may have several kinds of atoms,
but we now show that all of them originate the same kind of atomic pair:
\begin{prop}  \label{idatomic}
If $\langle \down e,\up e\rangle$ is an atomic pair of $\Sp X$, then $\down e$ and $\up e$ are retracts of representable functors
$\down x$ and $\up x$ for an object $x\in X$; furthermore, $\langle \down e,\up e\rangle$ can be obtained as the reflection of the 
atom ${\bf e}\tto^e X$ corresponding to an idempotent arrow $e:x\to x$ in $X$.
\end{prop}
\noindent So the choice of considering, in the definition of $\oX$ and $\Xo$, not all atomic pairs but only those which are reflections of actual atoms of $X$
was not restrictive, at least in this case.
\pf Note first that, by the above remarks, given $[\langle l,r\rangle]\in\ten(P,Q)$, $\alpha:P\to P'$ and $\beta:Q\to Q'$, 
we have $\ten(\alpha,\beta)[\langle l,r\rangle] = [\langle\alpha l,\beta r\rangle]\in\ten(P',Q')$. 
Now, let $u\in\ten(\down e,\up e) = \comp(\down e\,\times\!\up e)$ be biuniversal, and take an 
element $\langle lu, ru\rangle\in\,\down e\,\times\!\up e$ in this component: $[\langle lu,ru\rangle] = u$. 
Suppose furthermore that $\langle lu, ru\rangle$ is over $x\in X$ (that is, $lu\in(\down e)x$ and $ru\in(\up e)x$ ); by the (right) universality 
of $u$, there is a (unique) $\beta:\,\up e\to\up x$ such that $\ten(\down e,\beta)u = [\langle lu, {\rm id}_x\rangle]\in\ten(\down e,\up x)$; by Yoneda, 
there is a (unique) $\gamma:\,\up x\to\up e$ such that $\gamma:{\rm id}_x\mapsto ru\in(\up e)x$. So 
\begin{eqnarray*}
\ten(\down e,\gamma\circ\beta)u = \ten(\down e,\gamma)(\ten(\down e,\beta)[\langle lu,ru\rangle]) \\
= \ten(\down e,\gamma)[\langle lu,{\rm id}_x\rangle] = [\langle lu,ru\rangle] = u 
\end{eqnarray*}
but also $\ten(\down e,\up e)u = u$, and so, again by the (right) universality of $u$, $\gamma\circ\beta = \,\up e$ (the identity of $\up e$), 
as required. Of course, one proves symmetrically that $\down e$ is a retract of the representable $\down x$.
For the second part, observe that $\beta\circ\gamma : \,\up x\to \,\up x$ corresponds by Yoneda to an idempotent $e':x\to x$ in $X$,
whose reflection $\up e'$ is isomorphic to $\up e$, since both of them are retracts of $\up x$ associated to the same idempotent
in $[X,\Set]$ (see Proposition~\ref{k} above). Since an atomic pair is also a Dedekind cut (see Proposition~\ref{ded}), any of its components 
determines the other one (up to isomorphisms). So $\up e'\cong\,\up e$ implies $\down e'\cong\,\down e$ and the proof is complete.
\epf
Since a possible construction of the Cauchy completion of a category $X$ is given by the full subcategory of $[X\op,\Set]$
generated by the retracts of representable functors, and since Cauchy completion is defined up to equivalence 
(see e.g.~\cite{bor94} and~\cite{law89}), we have the following 
\begin{corol}  \label{cauchy}
In the space $\Sp X$, associated to the category $X$, $\,^{\rm o}(\Sp X)$ and $(\Sp X)^{\rm o}$ are the Cauchy completions of $X$ and $X\op$ respectively.
\end{corol}
\epf
Now we can return to one of the examples of Section~\ref{eval}: what is the evaluation of a closed part at an idempotent atom?
Recall that, as seen in Proposition~\ref{idat}, given a dof $D$ on $X$, the evaluation of $D$ at $e$ 
\[ \ev_e D := \hom(e,D) \cong \ten(e,D) \cong \hom(\up e,D)\cong \ten(\down e,D) \]
is the subset $D'x\subseteq Dx$ of elements fixed by the mapping $De:Dx\to Dx$, and similarly for open parts.
Then, we have the following generalizations of the Yoneda and co-Yoneda lemma (see Propositions~\ref{atomic3} and equation~\ref{coyo}):
\begin{eqnarray*} 
\hom(\down e,A) \cong A'x \qq && \qq \hom(\up e,D) \cong D'x  \\
\ten(A,\up e) \cong A'x  \qq && \qq \ten(\down e,D) \cong D'x
\end{eqnarray*}
which can be rephrased, for closed parts, as follows:
\begin{itemize}
\item
given an idempotent $e:x\to x$ in $X$ and a functor $D:X\to\Set$, the natural 
transformations $\up e\to D$ are in bijective correspondence with the elements of $Dx$ fixed by $e$;
\item
given an idempotent $e:x\to x$ in $X$ and a functor $D:X\to\Set$, the tensor product of functors $\down e\,\otimes D$
is given by the elements of $Dx$ fixed by $e$.
\end{itemize}
Both of them can be easily proved directly: in the former case, $\alpha:\,\up e\to D$ corresponds to $\alpha e\in D'x$, with inverse given
by $a\mapsto (D-)a$; in the latter $[\langle y\tto^f x,b\rangle]\in \down e\,\otimes D$ corresponds to $(Df)b\in D'x$, with 
inverse $a\mapsto [\langle x\tto^e x,a\rangle]$. 

In particular, given two atoms $e,e'\in\At X$ associated to idempotent arrows $e:x\to x$ and $e':x'\to x'$,
we have that $\oX(e,e') = \hom(\down e,\down e') \cong \ten(\down e',\up e)$, is the set of arrows $f:x\to x'\in(\down e')x$ fixed by $e$;
but by Proposition~\ref {k}, $f:x\to x'$ is in $(\down e')x$ if it is fixed by $e'$, so that $\oX(e,e')$ is the set of arrows $f:x\to x'$ 
such that $f\circ e = f = e'\circ f$. 
Thus we find the well-known description of the Cauchy completion of a category $X$ as its ``Karoubi envelope" (see e.g.~\cite{law89}),
which can be seen as the full subcategory of $\oX$ generated by the idempotent atoms, and that (by Propositions~\ref{idatomic})
is equivalent to $\oX$ itself.  
Furthermore, the interpretation of a df $A$ as a presheaf on $\oX$ (or on the Karoubi envelope of $X$) gives the (unique, up to isomorphisms) 
extension of the corresponding presheaf $A:X\op\to\Set$, whose value at any idempotent atom $e$, associated to $e:x\to x$, is the corresponding
fixed set $A'x\subseteq Ax$ (and similarly for the dof's). 

Finally, observe that the atomic pairs of $\Sp X$ coincide with the adjoint pairs of modules $L\adj R:\1\to X$ in the original definition
of the Cauchy completion of $X$ (see~\cite{law73} and~\cite{bor94}):
\begin{corol} \label{adjuni}
The modules $L:\1\to X$ and $R:X\to\1$ are adjoint iff there is a biuniversal element in $R\otimes L$ for the 
bifunctor $\otimes:\Mod(X,\1)\times\Mod(\1,X)\to\Set$.
\end{corol}
\epf
\noindent The unit $\eta:\, 1\to\,\down e\,\,\otimes\!\up e = \ten(\down e,\up e)$ of the 
adjunction $\down e\adj\,\,\up e:\1\to X$ selects the biuniversal element $e:x\to x$ among the endomorphisms $f:x\to x$ such that $f\circ e = f = e\circ f$.
The counit $\eps_{y,z}:(\up e\,\otimes\!\down e)(y,z)\to X(y,z)$ is given by composition. 

It seems likely that Corollary~\ref{adjuni} holds in the generic $\cal V$-valued context, but the details of the proof are still to be worked out.

\section{Reflection and coreflection in discrete fibrations}
\label{reflections}

In this section we complete the proof that $\Sp X = \langle \Cat/X,\open X,\closed X\rangle$ is a bipolar space,
showing directly that the formulas in Proposition~\ref{main} give indeed the reflection and the coreflection
of $\Cat/X$ in the full subcategory $\closed X$ of discrete opfibrations.
Of course, a dual proof works for discrete fibrations.

\subsection{The reflection in discrete opfibrations}

Like elsewhere, we shall follow the common ``abuse" of denoting a category $P\tto^{\pi}X$ over $X$ simply with $P$;
so, depending on the context, $P$ will indicate the functor $\pi$ or the corresponding total category.

We want to show that, for any discrete opfibration $D$ and $P\tto^{\pi}X$ in $\Cat/X$, there is a bijection 
\begin{equation} \label{eq1}
\frac{\up P \to D}{P\to D}
\end{equation}
where $\up P$ is given by Proposition~\ref{main} and Corollary~\ref{maincor}, that is it is defined on objects $x\in X$ by
\begin{equation}  \label{eq2}
(\up P)x = \ten(\down x,P) = \comp (\down x \times P) = \comp (P/x)
\end{equation}
Furthermore, such a bijection should be natural in $D\in\closed X$ and $P\tto^{\pi}X \in \Cat/X$.

\begin{prop}  \label{alpha1}
The morphisms over the line correspond (naturally) to the natural transformations in $[X,\Cat]$ between the functors 
\[ P/-:x\mapsto P/x \qq\qq D:x\mapsto\disc Dx \]
that is to the families of functors $\alpha_x:P/x\to \disc Dx$ such that
\[ h\alpha_x(a,f) = \alpha_y(a,h\circ f) \] 
for any $h:x\to y$ and $f:\pi a\to x$ in $X$.
\end{prop}
\pf Recall first that $\closed X$ is equivalent to $[X,\Set]$, and so the morphisms over the line can be taken as the natural transformations
\[ \alpha_x:(\up P)x\to Dx \] 
and each $\alpha_x$ corresponds by~(\ref{eq2}) and the adjuction $\comp\adj\disc$ to a functor 
\[ \alpha_x:P/x\to \disc Dx \] 
with the discrete category on $Dx$ as codomain; and since $[X,\Set]$ is reflective in $[X,\Cat]$, 
and the reflection being given pointwise by $\comp:\Cat\to\Set$, the naturality of the former family of mappings $\alpha_x$ 
coincides with the naturality of the corresponding family of functors.
\epf
On the other hand, the morphisms under the line are in $\Cat/X$, and so they should be in particular functors $\phi$ between 
the corresponding total categories; but since $D$ is an opfibration, these functors are in fact determined by their object mappings: 
\begin{prop}  \label{phi1}
A morphism $\phi:P\to D$ over $X$ is a family of mappings $\phi_x:\point Px\to Dx$ which associate an object-element in $Dx$
to any object in $P$ over $x$ in such a way that, for any $u:a\to b$ over $f:x\to y$
\[ \phi_y b = f(\phi_x a) \]
\end{prop}
\pf The condition is clearly necessary and one easily sees that it is also sufficient, that is, such a family of mappings defines
a functor over $X$.
\epf
To estabilish the desired natural bijection 
\begin{equation} \label{bijref}
\frac{P/- \to \disc D}{P\to D}
\end{equation} 
(which by Proposition~\ref{alpha1} is equivalent to~(\ref{eq1}))
we define, for any morphism $\phi:P\to D$ in $\Cat/X$, a family of mappings $\alpha(\phi)_x:\point(P/x)\to Dx$ as follows: 
\begin{equation} \label{bij1}
\alpha(\phi)_x(a,f) = f(\phi_{\pi a}a)
\end{equation} 
for any $a\in P$ and $f:\pi a\to x$.
In the other direction we define, for any morphism $\alpha:P/- \to \disc D$ in $[X,\Cat]$, a family of mappings $\phi(\alpha)_x:\point Px\to Dx$
as follows: 
\begin{equation}  \label{bij2}
\phi(\alpha)_x a = \alpha_x(a,{\rm id}_x)
\end{equation}
for any $a\in P$ over $x$.

\begin{prop}
Equations~{\rm(\ref{bij1})} and~{\rm(\ref{bij2})} define a natural bijection as required in equation~{\rm(\ref{bijref})}, 
and so $\ten(\down -,P)$ is indeed the reflection of $P\tto^{\pi}X$ in the opfibrations.
\end{prop}
\pf We must check that 
\begin{enumerate}
\item
the $\alpha(\phi)_x$ just defined give the components of a natural transformation $\alpha(\phi):P/- \to D$ in $[X,\Cat]$;
\item
the $\phi(\alpha)_x$ just defined are the ``components" of a morphism $\phi:P\to D$ in $\Cat/X$;
\item
the correspondences $\phi\mapsto\alpha(\phi)$ and $\alpha\mapsto\phi(\alpha)$ are each other's inverse:
\[ \phi = \phi(\alpha(\phi)) \qq\qq \alpha = \alpha(\phi(\alpha)) \]
\item
the bijection is natural.
\end{enumerate}

\begin{enumerate}
\item
First, the mappings $\alpha(\phi)_x:\point (P/x)\to Dx$ are in fact the object mappings of functors $\alpha(\phi)_x:P/x\to \disc Dx$
because, if $u:a\to b$ is an arrow $(a,f)\to (b,g)$ in $P/x$, then 
\begin{eqnarray*}
\alpha(\phi)_x(a,f) = \alpha(\phi)_x(a,g\circ\pi u) = (g\circ\pi u) (\phi_x a) = g((\pi u)(\phi_x a)) \\ 
= g(\phi_{\pi b} b) = \alpha(\phi)_x(b,g) 
\end{eqnarray*}
Furthermore, they are the components of a natural transformation $\alpha(\phi):P/- \to D$ in $[X,\Cat]$ 
because, for any $h:x\to y$ and $f:\pi a\to x$ in $X$, 
\[ h(\alpha(\phi)_x(a,f)) = h(f(\phi_x a)) = (h\circ f)(\phi_x a) = \alpha(\phi)_y(a,h\circ f) \]
 
\item
In the other direction, to show that $\phi(\alpha)_x:\point Px\to Dx$ defines a morphism $\phi:P\to D$ over $X$,
observe first that if $u:a\to b$ is an arrow in $P$, then $u$ itself is an arrow $(a,\pi u)\to (b,{\rm id}_{\pi b})$ in $P/x$,
and so the objects $(a,\pi u)$ and $(b,{\rm id}_{\pi b})$ must have the same image under any functor toward a discrete category.
By the naturality of $\alpha$ as in Proposition~\ref{alpha1} and by the above remark we then have
\[ \pi u (\phi(\alpha)_{\pi a} a) = \pi u (\alpha_{\pi a}(a,{\rm id}_{\pi a})) = \alpha_{\pi b}(a,\pi u) = 
\alpha_{\pi b}(b,{\rm id}_{\pi b}) = \phi(\alpha)_{\pi b} b \]
for any $u:a\to b$ in $P$, as required by Proposition~\ref{phi1}. 
\item
The correspondences $\phi\mapsto\alpha(\phi)$ and $\alpha\mapsto\phi(\alpha)$ are each other's inverse because, for any $a\in P$ over $x$,
\[ \phi(\alpha(\phi))_xa = \alpha(\phi)_x(a,{\rm id}_a) = {\rm id}_a(\phi_xa) = \phi_xa \]
and in the other direction, for any $f:\pi a\to x$ in $P$,
\[ \alpha(\phi(\alpha))_x(a,f) = f(\phi(\alpha)_{\pi a}a) = f(\alpha_{\pi a}(a,{\rm id}_a)) = \alpha_x(a,f) \]
again by the naturality of $\alpha$.
\item
Finally, it is straightforward to check that the bijection is indeed natural, and so the proposition is proved.
\end{enumerate}
\epf

\subsection{The coreflection in discrete opfibrations}
\label{coreflections}

The proofs in this subsection follow closely those of the previous one, displaying ``almost dual" aspects that we highlight at the end of the section.

We want to show that, for any discrete opfibration $D$ and $P\tto^{\pi}X$ in $\Cat/X$, there is a bijection 
\begin{equation}    \label{bijcoref}
\frac{D\to P\rup}{D\to P}
\end{equation}
where $P\rup$ is given by the formulas in Proposition~\ref{main}, that is it is defined on objects $x\in X$ by
\begin{displaymath}  
(P\rup)x = \hom(\up x,P) 
\end{displaymath}
Furthermore, such a bijection should be natural in $D\in\closed X$ and $P\tto^{\pi}X \in \Cat/X$.

The morphisms under the line are in $\Cat/X$, and so they should be in particular functors $\phi$ between 
the corresponding total categories; but since $D$ is an opfibration, we have: 
\begin{prop} \label{phi2}
A morphism $\phi:D\to P$ over $X$ is given by a mapping $\phi:\point D\to\point P$ over $X$, which associates an object $\phi a$ in $P$ over $x$ 
to any object-element $a\in Dx$, and by a mapping $\phi:(a,f)\mapsto \phi(a,f)$ which takes a pair $(a,f)$ with $a\in Dx$ and $f:x\to y$
in $X$ to an arrow $\phi(a,f)$ in $P$ over $f$ in such a way that 
\begin{enumerate}
\item
$\phi(a,f):\phi a\to \phi (fa)$
\item
$\phi(a,{\rm id}_{\pi a}) = {\rm id}_{\phi a}$ and $\phi(fa,g)\circ \phi(a,f) = \phi(a,g\circ f)$
\end{enumerate}
for any $a\in Dx$, $f:x\to y$ and $g:y\to z$.
\end{prop}
\pf Indeed, by the first condition, $\phi$ is a graph morphism $\phi:D\to P$ over $X$, while the second one states that it is in fact a functor.
\epf
\begin{prop} \label{alpha2}
The morphisms over the line in~{\rm(\ref{bijcoref})} correspond to the families of mappings 
\[ \alpha_x:Dx\to (P\rup)x \] 
which associate to each $a\in Dx$ a morphism $\up x\to P$ over $X$, that is a mapping
$\xi:f\mapsto\xi f$, which takes an arrow $f:x\to y$ in $X$ to an object $\xi f$ in $P$ over $y$, 
and a mapping $\xi:(f,g)\mapsto \xi(f,g)$ which takes a pair $(f,g)$ with $f:x\to y$ and $g:y\to z$
in $X$ to an arrow $\xi(f,g)$ in $P$ over $g$, in such a way that
\begin{enumerate}
\item
$\xi(f,g):\xi f\to \xi (g\circ f)$
\item
$\xi(f,{\rm id}_y) = {\rm id}_{\xi f}$ and $\xi(g\circ f,l)\circ \xi(f,g) = \xi(f,l\circ g)$
for any $f:x\to y$, $g:y\to z$ and $l:z\to w$.
\item
$(\alpha_y (ha))f = (\alpha_x a)(f\circ h)$ and $(\alpha_y (ha))(f,g) = (\alpha_x a)(f\circ h,g)$
for any $h:x\to y$, $f:y\to z$ and $g:z\to w$.
\end{enumerate}
\end{prop}
\pf Since $\closed X$ is equivalent to $[X,\Set]$, a morphism over the line can be taken as a natural transformation
$\alpha:D\to P\rup$, and since $\alpha_x:Dx\to (P\rup)x$ takes values in the morphisms $\up x\to P$ over $X$, the first
two conditions follow directly from Proposition~\ref{phi2}, with $D=\up x$.
Furthermore, as a presheaf $P\rup$ acts as follows: if $\xi\in{\hom}(\up x,P)$ and $h:x\to y$ in $X$, then $h\xi\in{\hom}(\up y,P)$
is defined by 
\[ (h\xi)f = \xi(f\circ h) \qq\qq (h\xi)(f,g) = \xi(f\circ h,g) \]
for any $f:y\to z$ and $g:z\to w$ in $X$.
(The reader may easily check that $h\xi$ is indeed still a morphism $\up y\to P$ over $X$.)
Then the naturality of the $\alpha_x:Dx\to P\rup$ is equivalent to the third condition of the proposition.
\epf
To estabilish the desired natural bijection~(\ref{bijcoref}) we define, for any morphism $D\to P$ over $X$ 
as in Proposition~\ref{phi2}, a family of mappings $\alpha(\phi)_x:Dx\to (P\rup)x$ as follows: 
\begin{equation} \label{bij3}
(\alpha(\phi)_x a)f = \phi(fa) \qq\qq   (\alpha(\phi)_x a)(f,g) = \phi(fa,g)
\end{equation} 
for any $a\in Dx$, $f:\pi a\to x$ and $g:y\to z$;
and in the other direction we define, for any morphism $\alpha:D\to P\rup$ in $[X,\Set]$ as in Proposition~\ref{alpha2}, 
the mapping $\phi(\alpha)$ as follows: 
\begin{equation}  \label{bij4}
\phi(\alpha) a = (\alpha_x a){\rm id}_x   \qq\qq  \phi(\alpha)(a,f) = (\alpha_x(a,f))({\rm id}_x,f)
\end{equation}
for any $a\in P$ over $x$ and $f:x\to y$.

\begin{prop}
Equations~{\rm(\ref{bij3})} and~{\rm(\ref{bij4})} define a natural bijection as required in equation~{\rm(\ref{bijcoref})}, 
and so ${\hom}(\up -,P)$ is indeed the coreflection of $P\tto^{\pi}X$ in the opfibrations over $X$.
\end{prop}
\pf We must check that 
\begin{enumerate}
\item
the $\alpha(\phi)_x$ just defined give the components of a natural transformation $\alpha(\phi):D\to P\rup$ in $[X,\Set]$;
\item
the $\phi(\alpha)$ just defined correspond to a morphism $\phi:D\to P$ in ${\Cat}/X$;
\item
the correspondences $\phi\mapsto\alpha(\phi)$ and $\alpha\mapsto\phi(\alpha)$ are each other's inverse:
\[ \phi = \phi(\alpha(\phi)) \qq\qq \alpha = \alpha(\phi(\alpha)) \]
\item
the bijection is natural.
\end{enumerate}

\begin{enumerate}
\item
First, the mappings $\alpha(\phi)_x$ have indeed values in $P\rup$, because the conditions of Proposition~\ref{phi2} on $\phi$ 
imply in a straightforward manner the corresponding ones on $\alpha(\phi)_x a$ required in Proposition~\ref{alpha2}.  
Furthermore, the naturality condition is fulfilled too, because for any $h:x\to y$, $f:y\to z$ and $g:z\to w$ in $X$ 
\[ (\alpha(\phi)_y(ha))f = \phi (f(ha)) = \phi((f\circ h)a) = (\alpha(\phi)_x a)(f\circ h) \]
\[ (\alpha(\phi)_y(ha))(f,g) = \phi (f(ha),g) = \phi ((f\circ h)a,g) = (\alpha(\phi)_x a)(f\circ h,g) \] 
\item
In the other direction, to show that $\phi(\alpha)$ defines a morphism $D\to P$ over $X$, we must check that the conditions
of Proposition~\ref{phi2} hold: 
\begin{enumerate}
\item 
for any $a\in Dx$ and $f:x\to y$, $\phi(\alpha)(a,f):\phi(\alpha) a\to \phi(\alpha) (fa)$.
That the domain of $\phi(\alpha)(a,f)$ is indeed $\phi(\alpha) a$ is immediate; as for the codomain, by the naturality of $\alpha$
\[ \phi(\alpha) (fa) = (\alpha_y (fa)){\rm id}_y = (\alpha_x a)f \]
and the latter is the codomain of $\phi(\alpha)(a,f) = (\alpha_x a)({\rm id}_x,f)$, by the corresponding property of $\alpha_x a$
in Proposition~\ref{alpha2}.
\item
that $\phi(\alpha)(a,{\rm id}_{\pi a}) = {\rm id}_{\phi(\alpha) a}$ is easily verified, while 
for any $a\in Dx$, $f:x\to y$ and $g:y\to z$
\begin{eqnarray*}
\phi(\alpha)(a,g\circ f) = (\alpha_x a)({\rm id}_x,g\circ f) = (\alpha_x a)(f\circ {\rm id}_x,g)\circ (\alpha_x a)({\rm id}_x,f) \\
= (\alpha_x (fa))({\rm id}_y,g)\circ (\alpha_x a)({\rm id}_x,f) = \phi(\alpha)(fa,g)\circ \phi(\alpha)(a,f)
\end{eqnarray*} 
again by the properties of $\alpha$ in Proposition~\ref{alpha2}.
\end{enumerate}

\item
The correspondences $\phi\mapsto\alpha(\phi)$ and $\alpha\mapsto\phi(\alpha)$ are each other's inverse because
for any $a\in Dx$, $f:x\to y$ and $g:y\to z$
\[ \phi(\alpha(\phi))a = (\alpha(\phi)_xa){\rm id}_a = \phi({\rm id}_a a) = \phi a \]
\[ \phi(\alpha(\phi))(a,f) = (\alpha(\phi)_xa)({\rm id}_a,f) = \phi({\rm id}_a a,f) = \phi (a,f) \]
and
\[ (\alpha(\phi(\alpha))_xa)f = \phi(\alpha)(fa) = (\alpha_x(fa)){\rm id}_y = (\alpha_x a)({\rm id}_y\circ f) = (\alpha_xa)f \]
\begin{eqnarray*}
(\alpha(\phi(\alpha))_xa)(f,g) = \phi(\alpha)(fa,g) = (\alpha_x(fa))({\rm id}_y,g) \\ 
= (\alpha_x a)({\rm id}_y\circ f,g) =(\alpha_xa)(f,g) 
\end{eqnarray*}
\item
Again, it is straightforward to check that the bijection is indeed natural, and so the proposition is proved.
\end{enumerate}
\epf
\noindent
If we consider only the object mappings of the morphisms $\phi$ and $\alpha a$ in the present section, 
the ``duality" of the proofs of the reflections and of the coreflections in discrete opfibrations is more evident, as summarized below.

\begin{itemize}
\item
Naturality conditions for $\alpha$: if $h:x\to y$
\begin{eqnarray*} 
h(\alpha_x(a,f)) &=& \alpha_y(a,h\circ f) \\ 
(\alpha_y (ha))f &=& (\alpha_x a)(f\circ h)
\end{eqnarray*}
\item
Definition of $\alpha(\phi)$ and $\phi(\alpha)$:
\begin{eqnarray*}  
\alpha(\phi)_x(a,f) = f(\phi_xa) \qq && \qq \phi(\alpha)_x a = \alpha_x(a,{\rm id}_a)  \\
(\alpha(\phi)_x a)f = \phi(fa) \qq && \qq  \phi(\alpha) a = (\alpha_xa){\rm id}_x 
\end{eqnarray*}
\item
Proof of the naturality of the $\alpha(\phi)$: if $h:x\to y$
\begin{eqnarray*}  
h(\alpha(\phi)_x(a,f)) = h(f(\phi_x a)) &=& (h\circ f)(\phi_x a) = \alpha(\phi)_y(a,h\circ f)  \\
(\alpha(\phi)_y(ha))f = \phi (f(ha)) &=& \phi((f\circ h)a) = (\alpha(\phi)_x a)(f\circ h)
\end{eqnarray*}

\end{itemize}

\section{The category of bipolar spaces}
\label{doctrines}

In this section we define the category $\Bip$ of bipolar spaces and continuous maps and consider various ``space functors" valued in it,
notably the (pseudo)functors $\Sp:\Cat\to\Bip$ and $\Sp:\Grf\to\Bip$; in the other direction, we define a ``base functor" $\Ba:\Bip\to\Cat$, 
by showing that atoms are preserved by (the left adjoint component of) a continuous map.
Then, connections with the Kan extensions of set functors are analyzed.

\subsection{Continuous maps}
\label{continuous}

The definition of continuous map is modelled on the effect that a functor $f:X\to Y$ has on the corresponding spaces $\Sp X$ and $\Sp Y$.
Recall that any such a functor gives rise to a pair of adjoint functors 
\[ f_!\adj f^*:\Cat/Y\to\Cat/X \] 
where $f_!$ is induced by composition with $f$, while $f^*$ is given by the pullback in $\Cat$.
The pair $f_!\adj f^*$ satisfies the Frobenius law, that is the (natural) morphism
\begin{equation} \label{frob}
\Phi_{P,Q}=f_!\pi_1\wedge(\eps_Q\circ f_!\p2) : f_!(P\times f^*Q) \to f_! P\times Q 
\end{equation}
is an isomorphism for any parts $P\in\Cat/X$ and $Q\in\Cat/Y$, $\eps$ being the counit of the adjunction. 
This follows essentially from the fact that the composition of two pullback squares, expressed by $f_!(P\times f^*Q)$, gives another
``outer" pullback square, expressed by $f_! P\times Q$. 
Furthermore, $P\in\Cat/X$ and $f_!P\in\Cat/Y$ have the same total category, so that
the functor $\comp^X:\Cat/X\to\Set$ factorizes through $f_!:\Cat/X\to\Cat/Y$ and $\comp^Y:\Cat/Y\to\Set$ up to the isomorphism
\[ \Lambda^f:\comp^X\iso\comp^Y\circ f_! \]
which takes the component $[a]_X$ of an object $a\in P$ to the component of the same object $[a]_Y$ in $f_!P$;
so, also $\disc_X$ factorizes through $\disc_Y$ and $f^*$, up to isomorphisms.
Finally, $f^*$ preserves open and closed parts:
if $A\in\Cat/Y$ is a df or a dof, so is $f^*A\in\Cat/X$.

So, if $X$ and $Y$ are bipolar spaces, we define a {\bf continuous map} $f:X\to Y$ as a pair of adjoint functors
\[ f_!\adj f^*:\P(Y)\to\P(X) \] 
which satisfies the Frobenius law, with $f^*$ preserving open and closed parts and with an isomorphism 
\[ \Lambda^f:\comp^X\iso\comp^Y\circ f_!:\P(X)\to\Set \]   

If $f:X\to Y$ and $g:Y\to Z$ are continuous maps, then $g\circ f:X\to Z$ is defined by composition of adjoint functors:
\[ g_!\circ f_!\adj f^*\circ g^*:\P(Z)\to\P(X) \] 
and it is again continuous; indeed, it satisfies the Frobenius law because, for any $P\in\P(X)$ and $R\in\P(Z)$,
the Frobenius map 
\[ \Phi_{P,R}: (g\circ f)_!(P\times (g\circ f)^*R) \to (g\circ f)_! P\times R \] 
factorizes through the isomorphisms $g_!(\Phi^f_{P,g^*R})$ and $\Phi^g_{f_!P,R}\,$:
\[ g_!f_!(P\times f^*g^*R) \iso g_!(f_!P\times g^*R) \iso g_!f_!P\times R \]
Furthermore, 
\[ \Lambda^{g\circ f}:\comp^X\iso\comp^Z\circ (g\circ f)_!:\P(X)\to\Set \] 
is given by the composition
\[ \comp^X\ttto^{\Lambda^f}\comp^Y\circ f_!\ttto^{\Lambda^g\circ f_!}\comp^Z\circ g_!\circ f_! \] 
and $(g\circ f)^*:\P(Z)\to\P(X)$ clearly preserves open and closed parts.
Thus we have the category $\Bip$, with bipolar spaces as objects and continuous maps as arrows.

By the remarks at the beginning of this section, there is a ``space functor", actually, a pseudofunctor
\[ \Sp:\Cat\to\Bip \]
and likewise, in the two-valued context, there is a space functor
\[ \Sp:\Pos\to\Bip \]
which takes a poset $X$ to its Alexandrov space, and a morphism $f:X\to Y$ of posets to the corresponding continuous map of
topological spaces $\Sp f : \Sp X \to \Sp Y$.

\subsection{The space functor for graphs}
\label{space}

We now define space (pseudo)functors
\[ \Sp:\Grf\to\Bip \qq\qq \Sp:\Grfo\to\Bip \]
where $\Grfo$ is the category of reflexive graphs, beginning with the object mapping of the former. So we present the bipolar space 
\[ \Sp X = \langle \Grf/X,\open X,\closed X\rangle \] 
associated to a graph $X$, straightforwardly generalizing the considerations of Section~\ref{graphs}.
\begin{enumerate}
\item
The category $\Grf/X$ has tensor: if $P\tto^{\pi}X$ is a graph over $X$, its points are the sections $X\to P$ of $\pi$, while its
components are those of the total graph $P$ (see Proposition~\ref{slice}).
\item
The full subcategory $\D_X$ of $\Grf/X$, generated by the nodes $x:{\rm D}\to X$ and the arrows $f:{\rm A}\to X$ of $X$,
is adequate for $\Grf/X$. In fact, $\Grf/X \cong [\D_X\op,\Set]$. (In general, if $\C = [\A\op,\Set]$ is a presheaf category,
one may consider an object $X\in\C$ as a df on $\A$, and so an object $Y\in\C/X$, being a df on $\A$ over the df $X$, is itself a df over the 
total category of $X$. But the latter is the category of elements ${\rm elts} X$ of the presheaf $X:\A\op\to\Set$; 
then the objects of $\C/X$ correspond to the df's on ${\rm elts} X$, that is $\C/X\cong [({\rm elts} X)\op,\Set]$ (see also~\cite{law89}). 
In the present case, $\D_X$ is ${\rm elts} X$, when the graph $X$ is considered as a presheaf on the category $\D$ of Section~\ref{endomappings}.)  
Multiplying $P\tto^{\pi}X$ by $x$ in $\Grf/X$ gives (as total graph) the fibre $Px$ over $x$, that is the set of nodes over $x$;
on the other hand, multiplying it by $f:x\to y$ gives a graph $Pf$ over the arrow $\rm A$, 
with $Pf(0)=Px$ and $Pf(1)=Py$ (see Section~\ref{mappings}).
So, an object $P\tto^{\pi}X$ of $\Grf/X$ is in particular a family of sets, the fibres $Px$ ($x\in X$), and a family of objects 
of $\Grf/{\rm A}$, the fibres $Pf$ ($f$ in $X$); and a morphism $\phi : P\to Q$ in $\Grf/X$ is determined by the corresponding
families of mappings between the fibres $\phi_x:Px\to Qx$ and of morphisms in $\Grf/{\rm A}$, $\phi_f:Pf\to Qf$ ($f:x\to y$ in $X$), 
with $\phi_f(0) = \phi_x$ and $\phi_f(1) = \phi_y$ (see also Section~\ref{catexp}).
\item
A part $P\in\P(X)=\Grf/X$ is closed or {\bf right functional} if it is orthogonal to any ``domain arrow" $x\to f$ in $\D_X$, that is if it is 
interpreted as a bijection by any domain arrow; similarly, $P$ is open or {\bf left functional} if it is interpreted as a bijection by any codomain arrow.
So $P$ is right (respectively left) functional iff all the $Pf$ above defined are right (respectively left) mappings as graph over $\rm A$; then in this case $P$
``is" a diagram (that is, a graph morphism) $X\to\Set$ (respectively $X\op\to\Set$), or equivalently a functor ${\hat X}\to\Set$ (respectively ${\hat X}\op\to\Set$)
from the free category generated by $X$. 
So, the full subcategories $\open X$ and $\closed X$ of open and closed parts may be seen as induced by the functors $\D_X\to {\hat X}$ which
reduce any domain or, respectively, codomain arrow in $\D_X$ to an identity. In particular, they are reflective and coreflective in $\Grf/X$.
\item 
As in the case of categories (see Section~\ref{categories}), the exponentials $P\imp Q$ in $\Grf/X$ can be computed as 
follows: $(P\imp Q)x = \Set(Px,Qx)$ and the arrows over $f$ in $P\imp Q$ are the morphisms $Pf\to Qf$ in $\Grf/{\rm A}$.  
In particular, if $P\in\open X$ and $Q\in\closed X$, the $Pf$ are left mappings and the $Qf$ are right mappings as in Section~\ref{mappings}, 
and so $P\imp Q$ is itself right functional (and dually $Q\imp P$ is left functional).
\item
The atoms of $\Sp X$ are the nodes ${\rm D}\tto^x X$, and these are actually strong atoms. The reflection $\up x$ of the atom $x$ is given
by $(\up x) y = {\hat X}(x,y)$, while if $f:y\to z$ is an arrow in $X$ then $(\up x) f$ is given, as a right mapping, 
by ${\hat X}(x,f):{\hat X}(x,y)\to{\hat X}(x,z)$. So, the categories $\,^{\rm o}(\Sp X)$ and $(\Sp X)^{\rm o}$ are (equivalent to) the 
free categories on $X$ and on $X\op$, and the associated interpretations take a left (respectively right) functional graph over $X$ to the 
corresponding presheaf on ${\hat X}$ (respectively, on ${\hat X}\op$). Thus $\Sp X$ is atomic, and we have the usual formulas for coreflections 
and reflections of a graph over $X$ in functional graphs, given by its interpretations and cointerpretations as a presheaf on the corresponding free category.
\end{enumerate}

Similarly, one defines the spaces $\Sp X = \langle \Grfo/X,\open X, \closed X \rangle$ associated to a reflexive graph $X\in\Grfo$. 
Note that now the atoms are the {\em points} of $X$, and that they are not strong.
At first sight reflexive graphs appear to be unsuitable for being analyzed by cofigures. Indeed, the very property that allows the 
construction of the homotopy category of any $\Grfo$-valued category (see~\cite{law86b} and~\cite{law89}), namely the fact that $\comp:\Grfo\to\Set$ 
preserves products, implies that the only information obtainable on a reflexive graph, by tensoring it with other graphs, is the number of its components;
this is in some sense analogous to the fact that the only information obtainable on an event, by knowing the probabilities of its conjunctions with other 
{\em independent} events, is its probability. 

In fact, as for categories, in the space $\Sp{\rm 1} = \langle \Grfo,\open{\rm 1}, \closed{\rm 1} \rangle$ 
the reciprocally coadequate subcategories $\open{\rm 1}$ and $\closed{\rm 1}$ of open and closed parts, reduce both to the subcategory
of discrete-constant graphs (so that $\Sp{\rm 1}$ is in fact a codiscrete Boolean space, as defined in Section~\ref{boole}; the same is true
for any codiscrete {\em category}). However, in general the spaces $\Sp X = \langle \Grfo/X,\open X, \closed X \rangle$
are as rich in content as those associated to irreflexive graphs: just consider the free reflexive graph on a given graph in $\Grf$.

Now we can turn to the definition of the space (pseudo)functor on graph morphisms $f:X\to Y$. As in the case of categories,
such a morphism gives rise to a pair of adjoint functors 
\[ f_!\adj f^*:\Cat/Y\to\Cat/X \] 
where $f_!$ is induced by composition, while $f^*$ is given by the pullback in $\Grf$.
The pair $f_!\adj f^*$ satisfies again the Frobenius law, and the functor $\comp^X:\Grf/X\to\Set$ factorizes through $f_!:\Grf/X\to\Grf/Y$ 
and $\comp^Y:\Grf/Y\to\Set$ up to isomorphisms:
\[ \Lambda^f:\comp^X\iso\comp^Y\circ f_! \]
Finally, $f^*$ preserves open and closed parts: if $A\in\Grf/Y$ is right or left functional, so is $f^*A\in\Grf/X$.
Analogous considerations hold for reflexive graphs.

Of course, also these space functors have a two-valued version:
\[ \Sp:\Rel\to\Bip \qq\qq \Sp:\Rel_0\to\Bip \]
where $\Rel$ and $\Rel_0$ are the categories of sets with an endorelation, or with a reflexive endorelation respectively.

\subsection{The base functor}
\label{base}

We now define a ``base" functor $\Ba:\Bip\to\Cat$ whose object mapping takes a space $X$ to the category $\oX$ defined in Section~\ref{atomic};
of course, one could consider $\Xo$ instead, and in fact a balanced choice would be to take the pair $\langle \oX , \Xo\rangle\in\Cat\times\Cat$. 
Recall that $\oX$ has the atoms in $\At X$ as objects, while $\oX(x,y) = \hom_X(\down x,\down y)$.
Then for any continuous map $f_!\adj f^*:\P(Y)\to\P(X)$ we want to define a functor 
\[ \Ba f : \oX \to \oY \]
and in particular its object mapping
\[ \Ba f : \At X \to \At Y \]
To do this, we first prove that the left adjoint $f_!:\P(X)\to\P(Y)$ takes atoms to atoms, so that we can define $(\Ba f)x = f_! x$, 
and then we extend it to a functor $\oX\to\oY$ as desired, thanks to the fact that $f^*$ preserves open parts. 
But before going into the details, we present the following result to illustrate the idea behind Proposition~\ref{bas}: 
\begin{prop}  \label{coadj2}
The left and right adjoint components of a continuous map are coadjoint functors, that is there are bijections
\[ \ten_X(P,f^*Q) \cong \ten_Y(f_!P,Q) \]
natural in $P\in\P(X)$ and $Q\in\P(Y)$.
\end{prop}
\pf We have 
\begin{eqnarray*}   
\ten_X(P,f^*Q) := \comp^X(P\times f^*Q) \cong \comp^Y f_!(P\times f^*Q) \\
\cong \comp^Y(f_!P\times Q) := \ten_Y(f_!P,Q) 
\end{eqnarray*}
where the first bijection is given by $\Lambda^f$ and the second by the Frobenius law, both of them natural.
\epf
Then if $x\in\P(X)$ is an atom, so is also $f_!x\in\P(Y)$, {\em roughly} because, for any $D\in\closed Y$ (and similarly for open parts),
\[ \ten_Y(f_!x,D) \cong \ten_X(x,f^*D) \cong \hom_X(x,f^*D) \cong \hom_X(f_!x,D) \]
where the first bijection is given by the above coadjunction, the second by the fact that $f^*D\in\closed X$ and $x$ is an atom,
and the third by the adjunction $f_!\adj f^*$.
But in fact we need to be more precise, exhibiting an actual biuniversal element $fu\in\ten_Y(f_!x,f_!x)$ associated to the biuniversal
element $u\in\ten_X(x,x)$, as required in the definition of atom.
The right definition of $fu$ is the most natural one:
\begin{equation} \label{fu}
fu := \comp^Y(f_!\pl\wedge f_!\p2)\Lambda^f_{x\times x} u 
\end{equation}
Observe that $fu$ is indeed in $\ten_Y(f_!x,f_!x)$, because $\Lambda^f_{x\times x}:\comp^X(x\times x) \to \comp^Y f_!(x\times x)$
and since $f_!\pl\wedge f_!\p2 : f_!(x\times x) \to f_!x\times f_!x$, 
\[ \comp^Y(f_!\pl\wedge f_!\p2) : \comp^Y f_!(x\times x) \to \comp^Y(f_!x\times f_!x) = \ten_Y(f_!x,f_!x) \]
To prove the biuniversality of $fu$ we must show that $\beta\mapsto\ten_Y(f_!x,\beta)fu$ is a bijection
\[ \ten_Y(f_!x,-)fu : \hom_Y(f_!x,D) \to \ten_Y(f_!x,D) \]
for any closed part $D$, and similarly for open parts. But
\begin{eqnarray*}  
\ten_Y(f_!x,\beta)fu = \comp^Y(f_!x\times\beta) \comp^Y(f_!\pl\wedge f_!\p2)\Lambda^f_{x\times x} u \\
= \comp^Y((f_!x\times\beta)\circ(f_!\pl\wedge f_!\p2))\Lambda^f_{x\times x} u = \comp^Y(f_!\pl\wedge(\beta\circ f_!\p2))\Lambda^f_{x\times x} u
\end{eqnarray*}
and so we have to show that
\begin{equation} \label{toprove}
\beta \mapsto \comp^Y(f_!\pl\wedge(\beta\circ f_!\p2))\Lambda^f_{x\times x} u 
\end{equation}
is a bijection. We need the following
\begin{lemma}  
The following equality of morphisms $f_!(x\times x) \to f_!x\times D$ in $\P(Y)$ holds: 
\[ f_!\pl\wedge(\beta\circ f_!\p2) = \Phi_{x,D}\circ f_!(x\times (f^*\beta\circ\eta_x)) \]
where $\Phi_{x,D}:f_!(x\times f^*D) \to f_! x\times D$ is the Frobenius map~{\rm(\ref{frob})} and $\eta$ is the unit of
the adjunction $f_!\adj f^*$. 
\end{lemma}
\pf Composing out both morphisms with the projections $\pl':f_!x\times D\to f_!x$ and $\p2':f_!x\times D\to D$ we get the same result:
in the first case we clearly obtain $f_!\pl$ and $\beta\circ f_!\p2$ respectively; on the other hand, denoting 
by $\pl'':x\times f^*D\to x$ and $\p2'':x\times f^*D\to f^*D$ the projections in the definition~(\ref{frob}) of $\Phi_{x,D}$, we have
\begin{eqnarray*} 
\pl'\circ\Phi_{x,D}\circ f_!(x\times (f^*\beta\circ\eta_x)) = \pl'\circ(f_!\pl''\wedge(\eps_D\circ f_!\p2''))\circ f_!(x\times (f^*\beta\circ\eta_x))  \\
= f_!\pl''\circ f_!(x\times (f^*\beta\circ\eta_x)) = f_!(\pl''\circ(x\times (f^*\beta\circ\eta_x))) = f_!\pl 
\end{eqnarray*}
\begin{eqnarray*} 
\p2'\circ\Phi_{x,D}\circ f_!(x\times (f^*\beta\circ\eta_x)) = \p2'\circ(f_!\pl''\wedge(\eps_D\circ f_!\p2''))\circ f_!(x\times (f^*\beta\circ\eta_x))  \\
= \eps_D\circ f_!\p2''\circ f_!(x\times (f^*\beta\circ\eta_x)) = \eps_D\circ f_!(\p2''\circ(x\times (f^*\beta\circ\eta_x))) \\ 
= \eps_D\circ f_!(f^*\beta\circ\eta_x\circ\p2) = \eps_D\circ f_!f^*\beta\circ f_!\eta_x\circ f_!\p2 \\
= \beta\circ\eps_{f_!x}\circ f_!\eta_x\circ f_!\p2 = \beta\circ f_!\p2
\end{eqnarray*}
where, in the last two passages, the naturality of the counit $\eps$ and the triangular identities have been used.
\epf 
We are now in a position to prove 
\begin{prop} \label{bas}
If $x$ is a (strong) atom of $X$, with biuniversal element $u\in\ten_X(x,x)$, and if $f:X\to Y$ is a continuous map, then $f_!x$ is
a (strong) atom of $Y$, with biuniversal element $fu\in\ten_Y(f_!x,f_!x)$, as defined in equation~{\rm(\ref{fu})} above.
\end{prop}
\pf To check that~(\ref{toprove}) is indeed a bijection, first apply the above lemma, then observe that the
morphism $f^*\beta\circ\eta_x$ is $\bar\beta:x\to f^*D$, corresponding to $\beta:f_!x\to D$ in the adjunction $f_!\adj f^*$, and use
the naturality of $\Lambda^f$ to get
\begin{eqnarray*} 
\comp^Y(f_!\pl\wedge(\beta\circ f_!\p2))\Lambda^f_{x\times x} u = \comp^Y(\Phi_{x,D}\circ f_!(x\times (f^*\beta\circ\eta_x)))\Lambda^f_{x\times x} u \\
= \comp^Y(\Phi_{x,D}\circ f_!(x\times\bar\beta))\Lambda^f_{x\times x} u = \comp^Y\Phi_{x,D}\circ\comp^Y f_!(x\times\bar\beta)\circ\Lambda^f_{x\times x} u \\
= \comp^Y\Phi_{x,D}\circ\Lambda^f_{x\times f^*D}\circ\ten_X(x\times\bar\beta) u
\end{eqnarray*}
Now the fact that
\[ \beta \mapsto \comp^Y\Phi_{x,D}\circ\Lambda^f_{x\times f^*D}\circ\ten_X(x\times\bar\beta) u \]
is a bijection $\hom_Y(f_!x,D) \to \ten_Y(f_!x,D)$ follows from its factorization through the bijections 
\[ \hom_Y(f_!x,D) \iso \hom_X(x,f^*D) \]
by the adjunction $f_!\adj f^*$,
\[ \hom_X(x,f^*D) \iso \ten_X(x,f^*D) \]
by the universality of $u\in\ten_X(x,x)$ and because $f^*D$ is closed in $X$,
\[ \ten_X(x,f^*D) = \comp^X(x\times f^*D)\iso \comp^Y f_!(x\times f^*D) \]
because $\Lambda^f$ is a natural isomorphism, and
\[ \comp^Y f_!(x\times f^*D) \iso \comp^Y (f_!x\times D) = \ten_Y(f_!x,D) \]
since the Frobenius map $\Phi_{x,D}$ is supposed to be an isomorphism, and then also $\comp^Y\Phi_{x,D}$ is such. 
Of course, the same proof works for open parts and for strong atoms too.
\epf
To define the functor $\Ba f : \oX \to \oY$ on arrows, we must define arrow mappings $\oX(x,y) \to \oY((\Ba f)x,(\Ba f)y)$, that is
\[ \Ba_{x,y}:\hom_X(\down x,\down y) \to \hom_Y(\down f_!x,\down f_!y)\]
These are obtained by composing 
\begin{eqnarray*} 
\hom_X(\down x,\down y) = \open X(\down x,\down y) \cong \int_{A\in\open X}(\open X(\down y,A),\open X(\down x,A)) \\
\cong \int_{A\in\open X}(\open X(y,A),\open X(x,A)) \to \int_{B\in\open Y}(\open X(y,f^*B),\open X(x,f^*B)) \\
\cong \int_{B\in\open Y}(\open Y(f_!y,B),\open Y(f_!x,B)) \cong \int_{B\in\open Y}(\open Y(\down f_!y,B),\open Y(\down f_!x,B)) \\
\cong \open Y(\down f_!x,\down f_!y) = \hom_Y(\down f_!x,\down f_!y)
\end{eqnarray*}
where the only non-bijective mapping is given by the composition of natural transformations with the functor $f^*:\open Y\to\open X$.
\begin{prop} \label{bas2}
If $f:X\to Y$ is a morphism of bipolar spaces, the mappings $\Ba_{x,y}:\oX(x,y)\to\oY((\Ba f)x,(\Ba f)y)$ do define a functor $\Ba f: \oX\to\oY$. 
Furthermore, 
\[ \Ba(\id_X) \cong \id_{\oX} \qq\qq \Ba(g\circ f) \cong \Ba g\circ\Ba f \]
for any $g:Y\to Z$ in $\Bip$.
\end{prop}
\pf In fact, $\Ba f$ is essentially the restriction to the atoms reflections of the functor $\down(-)\circ f_!:\open X\to\open Y$;
to see this, observe that for any $B\in\open Y$
\begin{eqnarray*}  
\hom_Y(\down(f_!\down x),B)\cong \hom_Y(f_!\down x,B)\cong \hom_X(\down x,f^*B) \\
\cong  \hom_x(x,f^*B) \cong \hom_Y(f_!x,B) \cong \hom_Y(\down(f_!x),B)
\end{eqnarray*}
and then $\down(f_!\down x)\cong\down f_!x$.
\epf
Recall that, by Corollary~\ref{cauchy}, given a category $X$, $\Ba(\Sp X)$ is the Cauchy completion of $X$. So, the (pseudo)functor
\[ \Ba\circ\Sp:\Cat\to\Cat \]
may be called the ``Cauchy functor". In the two-valued context, the Cauchy functor  
\[ \Ba\circ\Sp:\Pos\to\Pos \]
when restricted to partially ordered sets, is equivalent to the identity, since in this case $\Ba(\Sp X)$ is isomorphic to $X$ 
(while for general preordered sets we have only an equivalence $\Ba(\Sp X)\simeq X$). 
Finally, in the case of graphs, the (pseudo)functors 
\[ \Ba\circ\Sp:\Grf\to\Cat \qq\qq \Ba\circ\Sp:\Grfo\to\Cat \]
serve as ``free category functors".

\subsection{Kan extensions}
\label{kan}

If $A\in\Cat/Y$ is a df on $Y$ and $f:X\to Y$, then the df $f^*A\in\Cat/X$ on $X$, considered as a presheaf, is obtained by substitution
in the presheaf corresponding to $A$:
\[ (f^*A)x = A(fx) \qq\qq (f^*A)\alpha = A(f\alpha):(f^*A)y\to(f^*A)x \]
for any $\alpha:x\to y$ in $X$. Similarly for the dof's:
\[ (f^*D)x = D(fx) \qq\qq (f^*D)\alpha = D(f\alpha):(f^*D)x\to(f^*D)y \]
Since the Kan extensions of set functors refer to left and right adjoint to such a functor
\begin{equation}  \label{subst}
f^*:\closed Y\to\closed X
\end{equation}
obtained by restricting $f^*:\Cat/Y\to\Cat/X$ to the closed parts of $\Sp Y$, the formalism of bipolar spaces gives some indications on them too.

The case of right Kan extensions is more direct: if the functor~(\ref{subst}) has a right adjoint 
\[ f^*\adj\forall_f:\closed X\to\closed Y \]
then
\[ (\forall_fD)y \cong \hom_Y(y,\forall_fD) \cong \hom_Y(\up y,\forall_fD) = \closed Y(\up y,\forall_fD) \cong \closed X(f^*\up y,D) \]
for any $y\in Y$; since $(f^*\up y)x = (\up y)(fx) = Y(y,fx)$ we get 
\[ (\forall_fD)y \cong \hom_X(Y(y,f-),D) \cong \int_x Y(y,fx)\imp Dx \]
the end formula for the right Kan extension of $D$ along $f:X\to Y$.

More interesting is the case of left Kan extensions. 
Observe first that the functor~(\ref{subst}) factorizes as the inclusion $\iota:\closed Y\hookrightarrow\Cat/Y$ followed
by (the corestriction of) the pullback functor $f^*:\Cat/Y\to\Cat/X$, and so its left adjoint
\[ \exists_f\adj f^*:\closed Y\to\closed X \]
can be obtained as the composition of (the restriction of) $f_!\adj f^*$ followed by $\up(-)\adj\iota$:
\[ \exists_f = \up(-)\circ f_!:\closed X\to\closed Y \]
so that $\exists_fD = \up f_!D$; using the reflection formulas of Proposition~\ref{main}
and the coadjunction relation of Proposition~\ref{coadj2} we get 
\[ (\exists_fD)y = (\up f_!D)y \cong\ten_Y(\down y,f_!D) \cong\ten_X(f^*\!\down y,D) \] 
that is we have the coend formula for the left Kan extension of $D$ along $f:X\to Y$ (see Proposition~\ref{coend}):
\[ (\exists_fD)y \cong\ten_X(Y(f-,y),D) \cong \int^x Y(fx,y)\times Dx \]

\subsection{An example}
\label{ex}

With respect to the adjunction
\[ \hat{-} \adj |-|:\Cat\to\Grf \] 
the following relations between the spaces $\Sp X$ associated to categories or to graphs are easily seen to hold:
\begin{enumerate}
\item
For any $X\in\Grf$, the open and closed part of $\Sp X$ are the ``same" as those of $\Sp\hat X$.
\item
For any $X\in\Cat$, the open and closed part of $\Sp X$ are ``included" as full subcategories in those of $\Sp |X|$.
Furthermore, this inclusion preserves the tensor functor too. So, ultimately, the calculations in $\Sp X$ (such as those relative to Kan
extensions) can actually be carried out in $\Grf/|X|$. 
\end{enumerate}

For instance, let $\N$ be the monoid of natural numbers, $\Z$ the group of integers and $\Z_n$ the cyclic group of order $n$.
By the first remark, the open and closed part of $\Sp\N$ are the same as those of $\Sp\L$, that is the presheaves in $\open\N$
and in $\closed\N$ can be seen as left and right endomapping in $\Grf$. 
On the other and, $\Sp\Z$ is a Boolean space, whose clopen part are the same as those of $\Sp\N$ (see Section~\ref{bifibrations}),
and so also as those of $\Sp\L$. Furthermore the quotient homomorphism $f_n:\Z\to\Z_n$ induces a full and faithful functor $f_n^*:\clopen\Z_n\to\clopen\Z$,
so that ultimately $\clopen\Z$ and $\clopen\Z_n$ are included as full categories in $\Grf$, and the inclusion preserves also the tensor functor
(see Corollary~\ref{clopencor}).
A part $A\in\clopen\Z$ is easily visualized as an object of $\Grf$: it is a sum of cycles $A = \sum_{k=1}^\infty n_k\cdot \L_k$, where $n_k$ is the 
number (or set) of cycles $\L_k$ of lenght $k$; $A$ is in $\clopen\Z_n$ if the lenght of any of its cycles is a divisor of $n$.

Thus, the reflection and coreflection of $\clopen\Z$ in $\clopen\Z_n$, that is the left and right extensions along $f_n$, can be computed in $\Grf$
and give also the reflection and coreflection of the left-right endomappings in those of order $n$. 
Denoting with $*$ the object of $\Z_n$, observe that $f^*\up *$ and $f^*\!\down *$ are both the graph $\L_n$, to be thought of as a right and a left
endomapping respectively. The formulas of Section~\ref{kan} now give that $\forall_{f_n}A$ is the endomapping (of order $n$) on the 
cycles $\hom(\L_n,A)$ given by the (``clockwise") shift of $\L_n$; similarly $\exists_{f_n}A$ is the endomapping (of order $n$) on the 
``cocycles" $\ten(\L_n,A)$ of $A$, given by the (``counterclockwise") shift of $\L_n$.  

But $\hom(\L_n,-)$ and $\ten(\L_n,-)$ both preserve sums, the former because $\L_n$ is connected and the latter by Proposition~\ref{slice} 
(or also because it has the right adjoint $\neg\L_n$); thus, if $A\in\clopen\Z$ the calculus of $\hom(\L_n,A)$ and $\ten(\L_n,A)$ reduce to those 
of $\hom(\L_n,\L_k)$ and $\ten(\L_n,\L_k)$, and one easily see that $\hom(\L_n,\L_k) = k$ if $k$ is a divisor of $n$ and $0$ otherwise,
while $\L_n\times \L_k = \gcd(n,k)\cdot \L_{\lcm(n,k)}$, so that $\ten(\L_n,\L_k) = \gcd(n,k)$. Summarizing:
\begin{enumerate}
\item
The coreflection of bijective endomappings in those of order $n$ is obtained by taking only the cycles whose lenght is a divisor of $n$.
\item
The reflection of bijective endomappings in those of order $n$ is obtained by ``reducing" each $k$-cycle to a $\gcd(n,k)$-cycle. 
\end{enumerate}
Although the result may appear shallow when compared with the machinery employed, this example shows once more the conceptual advantage 
of an explicit consideration of the $\ten$ functor, along with the ``dual" $\hom$ functor.
Of course, similar ``visualizable" reasonings  can be applied to more complex situations as well.

\section{Appendix}

\subsection{Discrete bifibrations}
\label{bifibrations}

If $X$ is a topological space, the relation
\begin{equation}  \label{comp}
x\sim y \iff x\in B\imp y\in B,\q\forall B\in\clopen X 
\end{equation}
where $\clopen X$ is the set of clopen parts of $X$, is easily verified to be an equivalence relation on $X$, whose classes are the
connected components. If the points $x\in X$ have reflections $\clopen x$ in $\clopen X$ (that is, if the smallest clopen set $\clopen x$
containing $x$ exists) the~(\ref{comp}) becomes
\[ x\sim y \iff \clopen x\subseteq\clopen y \iff x\in\clopen y \]
In particular, for the Alexandrov space associated to a poset $X$ we get the grupoidal reflection of $X$ (see Section~\ref{preorders}).

We now briefly investigate those bipolar spaces for which parts have clopen reflections and coreflections
(see Proposition~\ref{clopen}), the key example being again the space $\Sp X$ associated to a category $X$.
Observe first that the negation $\neg A$ of a clopen part $A\in\clopen X$ of a bipolar space $X$ is itself clopen;
then the proof of Proposition~\ref{pcoadj} is still valid, and if we denote by $\bi(-):\P(X)\to\clopen X$ and $(-)\rbi:\P(X)\to\clopen X$ 
the reflection and the coreflection in clopen parts, we have (see Corollary~\ref{coadj} of Section~\ref{bip}):
\begin{corol} 
If $P$ and $Q$ are parts of a bipolar space, there are natural isomorphisms
\[ \ten(P,\bi Q) \cong \ten(\bi P,\bi Q) \cong \ten(\bi P,Q) \]
\[ \hom(P,Q\rbi) \cong \hom(\bi P,Q\rbi) \cong \hom(\bi P,Q) \]
\end{corol}
If $X$ has clopen reflections, we can define $\oXo$ the category with objects in $\At X$ and arrows in $\clopen X$, that is
\[ \oXo(x,y) = \hom(\bi x,\bi y) \]
and there is the obvious full and faithful functor $\lpr:\oXo\to\P(X)$.
Then argumentations similar to those of Section~\ref{atomic} show that $\oXo$ is self-dual via an isomorphism 
\[ \sigma:\oXo\to(\oXo)\op \]
with the identity as object mapping (so that in the two-valued case $\oXo$ is an equivalence relation), and that the interpretation and cointerpretation
of a clopen part of $X$ via $\lpr$ are the same, modulo this duality. If $X$ is atomic, $\oXo$ is adequate and coadequate for clopen
parts, and we have the formulas for reflection and coreflection of a part $P\in\P(X)$ in clopen parts (see Corollary~\ref{atomic5} of 
Section~\ref{atomic}):
\begin{eqnarray*} 
(\bi P) x \cong \ten(\bi x,P) \qq && \qq (\bi P) \alpha \cong \ten(\alpha,P) \\
(P\rbi) x \cong \hom(\bi x,P) \qq && \qq (P\rbi) \alpha \cong \hom(\alpha,P) 
\end{eqnarray*}
For the bipolar space $\Sp X$, $X\in\Cat$, we have that $\,^{\rm o}(\Sp X)^{\rm o}$ is equivalent to the grupoid reflection 
of $\,^{\rm o}(\Sp X)$ or $(\Sp X)^{\rm o}$; restricting to the objects-atoms, we get the grupoid reflection $G_X$ of $X$ or $X\op$,
and the above formulas give the reflection and coreflection of a category over $X$ in discrete bifibrations (that is,
in presheaves which act by bijections) as presheaves on $G_X$.  

\subsection{Comappings}
\label{comappings}

In our examples of bipolar spaces, open and closed parts were defined by orthogonality with respect to special morphisms in $\P(X)$.
In fact, in each $\P(X)$ we could find special (subterminal) objects, the ``nodes" and the ``arrows", such that any ``arrow" $f$ has ``domain" and
``codomain" morphisms $x\to f$ and $y\to f$ for uniquely determined ``nodes" $x$ and $y$; closed and open parts were then defined by orthogonality 
with respect to these domain and codomain morphisms. 
In particular, in Section~\ref{mappings} we defined right and left mappings of  by orthogonality with respect to $\delta_0:{\rm D_0}\to{\rm A}$ 
and $\delta_1:{\rm D_1}\to{\rm A}$ in $\Grf/{\rm A}$; e.g. $P = P\tto^{\pi}{\rm A}$ is a right mapping iff it is interpreted as a bijection by the 
arrow $\delta_0$. Since a recurrent topic of the present paper is to work with cofigures too, one is led to wonder what is obtained with the
corresponding ``co-orthogonality condition", namely that $P$ is cointerpreted as a bijection by the arrow $\delta_1:{\rm D_1}\to{\rm A}$.
The answer is: a (right) mapping with possibly repeated arrows; we call such an object of $\Grf/{\rm A}$ a (right) ``comapping".
Comappings do have the hall-mark of functionality: to any element $x$ in $P(0)$ there corresponds one and only one element $Px$ of the
codomain $P(1)$. The two interpretations of Section~\ref{mappings} still hold, but the first one now sounds 
``follow {\em any} of the arrows starting from $x$". 
On the other hand, comappings lack the unicity of mappings, since two of them may have the same ``functional role"; in particular, left and right
comappings cannot be reciprocally coadequate.

\subsection{Abstractions}
\label{abstractions}

The space functors $\Sp:\Cat\to\Bip$ (or, in a two-valued context, $\Sp:\Pos\to\Bip$), $\Sp:\Grf\to\Bip$ and $\Sp:\Grfo\to\Bip$ of Sections~\ref{continuous} 
and~\ref{space} are all instances of ``bipolar doctrines", i.e. pseudofunctors $\T\to\Bip$ for categories $\T$ with suitable properties (see~\cite{law70}).
Of course, these doctrines constitute a context which makes it possible generalizations of categorical concepts such as Kan extensions or 
Cauchy completions. Abstracting further, one can define temporal doctrines wherein, for each type $X\in\T$, the category of attributes or parts $\P(X)$ 
has four temporal operator, i.e. the endofunctors
\[ \down(-),\up(-),(-)\rdown\, ,(-)\rup\,:\P(X)\to\P(X) \] 
which fulfil the adjunction laws $\down(-)\adj (-)\rdown$ and $\up(-)\adj (-)\rup$ and the coadjunction law (see equations~(\ref{f1}) and~(\ref{f2})
in the Introduction).

Anyway, a general principle which seems to emerge from the present work is the following:
while in a two-valued context open and closed parts (generally speaking) are often disciplined by a strict duality (due to the self-duality 
of $\Two$ itself), and so one can dispense with the explicit consideration of both of them, this is not the case in a set-valued context:
the open and closed parts need to be dealt with simultaneously, their interaction being given not only by the fact that they belong to
the same category, but, most importantly, by the tensor functor. A point that may deserve to be deepened is whether a significant 
theory can be achieved by considering only open and closed parts, apart from their inclusion in an ampler category of parts, along with a 
bifunctor $\ten:\open X\times\closed X\to\Set$ satisfying suitable properties (such as a the existence of negation functors and reciprocal coadequacy).

\section{Conclusions}
\label{conclusions}

We hope to have shown that the conceptual frame of bipolar spaces offers a new perspective on important aspects of 
category theory, such as the Cauchy completion of a category, the tensor product (coend) of functors, the Kan extensions of set functors
and on the concept of mapping itself;
then it seems likely that a ``logic with components", wherein the tensor functor and the resulting notion of cofigure
play a major role, could be considered as a useful extension or deepening of ``generalized logic" (\cite{law73}).
Whether bipolar spaces have other theoretical aspects that deserve to be developed, or encompass more relevant examples
than those presented here (perhaps, valued in ``truth values" categories $\cal V$ other than $\Set$ or $\Two$), is open to further investigations.

\begin{refs}

\bibitem[Borceux, 1994]{bor94} F. Borceux, {\em Handbook of Categorical Algebra 1 (Basic Category Theory)}, 
Encyclopedia of Mathematics and its applications, vol. 50, Cambridge University Press, 1994.

\bibitem[Bunge \& Niefield, 2000]{bun00} M. Bunge and S. Niefeld, Exponentiability and Single Universes, 
{\em J. Pure Appl. Algebra} {\bf 148} (2000) 217-250.

\bibitem[Eilenberg \& Kelley, 1966]{eik66} S. Eilenberg and G.M. Kelly, {\em Closed Categories}, Proceedings of the Conference on Categorical Algebra, 
La Jolla 1965, Springer Verlag, 1966.
 
\bibitem[El Bashir \& Velebil, 2002]{elv02} R. El Bashir and J. Velebil, Simultaneously Reflective and Coreflective Subcategories of Presheaves, 
{\em Theory and Appl. Cat.} {\bf 10} (2002) 410-423.

\bibitem[Johnstone, 1999]{joh99} P. Johnstone, A Note on Discrete Conduch\'e Fibrations, 
{\em Theory and Appl. Cat.} {\bf 5} (1999) 1-11.
 
\bibitem[Lawvere, 1970]{law70} F.W. Lawvere, {\em Equality in Hyperdoctrines and the Comprehension Scheme as an Adjoint Functor},
Proceedings of the AMS Symposium on Pure Mathematics, XVII, 1970, 1-14. 

\bibitem[Lawvere, 1973]{law73} F.W. Lawvere, Metric Spaces, Generalized Logic and Closed Categories,
{\em Rend. Sem. Mat. Fis. Milano} {\bf 43} (1973) 135-166. Republished in {\em Reprints in Theory and Appl. Cat.} No. 1 (2002) 1-37.

\bibitem[Lawvere, 1986]{law86} F.W. Lawvere, Taking Categories seriously, 
{\em Revista Colombiana de Matematicas} {\bf 20} (1986) 147-178. Republished in {\em Reprints in Theory and Appl. Cat.} No. 8 (2005) 1-24.

\bibitem[Lawvere, 1986b]{law86b} F.W. Lawvere, Categories of Spaces may not be Generalized Spaces as exemplified by Directed Graphs, 
{\em Revista Colombiana de Matematicas} {\bf 20} (1986) 179-185. Republished in {\em Reprints in Theory and Appl. Cat.} No. 9 (2005) 1-7.

\bibitem[Lawvere, 1989]{law89} F.W. Lawvere, {\em Qualitative Distinctions between some Toposes of Generalized Graphs},
Proceedings of the AMS Symposium on Categories in Computer Science and Logic, Contemporary Mathematics, vol. 92, 1989, 261-299.

\bibitem[Lawvere, 1996]{law96} F.W. Lawvere, {\em Adjoints in and among Bicategories},
Proceedings of the 1994 Conference in Memory of Roberto Magari,
Logic \& Algebra, Lectures Notes in Pure and Applied Algebra, 180:181-189,
Ed. Ursini Aglian\`o, Marcel Dekker, Inc. Basel, New York, 1996. 

\bibitem[Mac Lane, 1971]{mac71} S. Mac Lane, {\em Categories for the Working Mathematician}, 
Graduate Texts in Mathematics, vol. 5, Springer, Berlin, 1971.

\bibitem[Mac Lane \& Moerdijk, 1991]{mam91} S. Mac Lane and I. Moerdijk, {\em Sheaves in Geometry and Logic (a First Introduction to Topos Theory)}, 
Universitext, Springer, Berlin, 1991.

\bibitem[Par\'e, 1973]{par73} R. Par\'e, Connected Components and Colimits,
{\em J. Pure Appl. Algebra} {\bf 3} (1973) 21-42.

\bibitem[Street, 2001]{str01} R. Street, Powerful Functors, expository note (2001), available at www.maths.mq.edu.au.

\end{refs}

\end{document}